\documentclass[a4paper,11pt]{article}

\usepackage{amssymb}
\usepackage{amsmath}
\usepackage{amscd}
\usepackage[english]{babel}
\usepackage[all]{xypic}
\usepackage[all]{xy}
\usepackage{url}
\usepackage{graphicx}
\usepackage{hyperref}

\usepackage{comment}

\newtheorem{theo}{Theorem}[section]
\newtheorem{prop}[theo]{Proposition}
\newtheorem{lem}[theo]{Lemma}
\newtheorem{cor}[theo]{Corollary}

\newtheorem{rema}[theo]{Remark}

\newcommand{\dem}[1] {\paragraph{Proof{#1}: }}

\def \Romannumeral #1 {\expandafter\uppercase\expandafter {\romannumeral #1} }

\def \fraco {{\rm{Frac\,}}}

\def \gal {{\rm{Gal}}}

\def \calo {{\cal O}}
\def \spec {{\rm{Spec\,}}}

\def \Hom {{\rm {Hom}}}

\def \ext {{\rm {Ext}}}

\def \id {{\rm{id\,}}}

\def \Z {{\bf Z}}
\def \Q {{\bf Q}}
\def \L {{\bf L}}
\def \F {{\bf F}}

\def \id {{\rm{id}}}
\def \G {{\bf G}_m}

\def\ov{\overline}

\def\smallsquare{\vbox{\hrule\hbox{\vrule height 1 ex\kern 1 ex\vrule}\hrule}}
\def\enddem{\hfill\smallsquare\vskip 3mm}
\def \merci {\paragraph{Acknowledgements. }}
\def \abstract{\paragraph{R\'esum\'e. }}
\def \abstractbis{\paragraph{Abstract. }}
\def \codams{\paragraph{AMS codes~: }}

\def \keywbis{\paragraph{Keywords~: }}

\DeclareFontFamily{U}{wncy}{}
\DeclareFontShape{U}{wncy}{m}{n}{%
   <5>wncyr5%
   <6>wncyr6%
   <7>wncyr7%
   <8>wncyr8%
   <9>wncyr9%
   <10>wncyr10%
   <11>wncyr10%
   <12>wncyr6%
   <14>wncyr7%
   <17>wncyr8%
   <20>wncyr10%
   <25>wncyr10}{}
\DeclareMathAlphabet{\cyrille}{U}{wncy}{m}{n}

\def \R{{\bf R}}

\def \et {{\textup{\'et}}}

\def\Cone{{\rm Cone}}

\def\Gm{{\mathbf{G}_m}}

\author{Cyril Demarche and David Harari}

\title{Artin--Mazur--Milne duality for fppf cohomology}

\begin{document}
\maketitle

\abstractbis We provide a complete proof of a duality theorem 
for the fppf cohomology of either a curve over a finite field 
or a ring of integers of a number field, which extends 
the classical Artin--Verdier Theorem in \'etale cohomology.
We also prove some finiteness and vanishing statements.

\codams 11G20, 14H25.

\keywbis fppf cohomology, arithmetic duality, Artin approximation 
theorem.

\section{Introduction}

Let $K$ be a number field or 
the function field of a smooth, projective, geometrically 
integral curve $X$ over a finite field. In the number field 
case, set $X:=\spec \calo_K$, where $\calo_K$ is the ring of integers 
of $K$. Let $U$ be a non empty Zariski open subset of $X$ and denote by
$N$ a commutative, 
finite and flat group scheme over $U$ with Cartier dual $N^D$.
Assume that the order of $N$ is invertible on $U$ (in particular
$N$ is \'etale).
The classical ``\'etale'' {\it Artin--Verdier Theorem} (cf.~\cite{MilADT}, Corollary 
II.3.3.) is a duality statement between 
\'etale cohomology $H^{\bullet}_{\et}(U,N)$ 
and \'etale cohomology with 
compact support $H^{\bullet}_{\et,c}(U,N^D)$. 
It has been known for a long time that this theorem is especially useful 
in view of concrete arithmetic applications: for example it yields 
a very nice method to prove deep results like Cassels--Tate duality for 
abelian varieties and schemes (\cite{MilADT}, section II.5) and their 
generalizations to $1$-motives (\cite{dhsza}, section 4); 
Artin--Verdier's Theorem 
also provides a ``canonical'' path to prove the Poitou--Tate's Theorem and its 
extension to complex of tori (\cite{demimrn}), which in turn turns out 
to be very fruitful to deal with local-global questions for (non 
necessarily commutative) linear algebraic groups (\cite{demlms}). 

\smallskip

It is of course natural to try to remove the condition that the order 
of $N$ is invertible on $U$. A good framework to do this is 
provided by fppf cohomology of finite and flat commutative 
group schemes over $U$, 
as introduced by J.S.~Milne in the third part of his book \cite{MilADT}. 
This includes the case of group schemes of order divisible by 
$p:={\rm Char} \, K$ in the function field case. 

Such a fppf duality theorem has been first announced by 
B.~Mazur\footnote{Thanks to A.~Schmidt for having pointed this out to us.}
(\cite{mazur}, Prop.~7.2), relying on work by M.~Artin and himself. 
Special cases have also been proved by M.~Artin and Milne 
(\cite{artmilne}). The precise statement 
of the theorem is as follows (see \cite{MilADT}, 
Corollary III.3.2.~for the number field case and Theorem III.8.2 for
the function field case):

\begin{theo} \label{thm AV}
Let $j : U \hookrightarrow X$ be a non empty open subscheme of $X$. Let 
$N$ be a finite flat commutative group scheme over $U$ with 
Cartier dual $N^D$.
For all integers $r$ with $0 \leq r \leq 3$, the canonical pairing
$$H^{3-r}_c(U, N) \times H^r(U, N^D) \to H^3_c(U, \G) \cong \Q / \Z$$
(where $H^r(U, N^D)$ is a fppf cohomology group and 
$H^{3-r}_c(U, N)$ a fppf cohomology group 
with compact support) induces a perfect duality between
the profinite group $H^{3-r}_c(U, N)$ and the 
discrete group $H^{r}(U, N^D)$.
Besides, these groups are finite in the number field case,
and they are trivial for $r \geq 4$ and $r < 0$ (resp.~for $r=3$ if $U \neq X$) 
in the function field case.
\end{theo}

For example, this extension of the \'etale Artin--Verdier Theorem is needed
to prove the Poitou--Tate exact sequence over global fields of characteristic 
$p$ (\cite{gonzcrelle}, Th.~4.8. and 4.11) 
as well as the Poitou--Tate Theorem over a global 
field without restriction on the order (\cite{cesnamrl}, Th.~5.1, which 
in turn is used in \cite{rosen}, \S 5.6 and 5.7). Results 
of \cite{MilADT}, section III.9.~(which rely on the fppf duality 
Theorem) are also 
a key ingredient in the proof of some cases
of the Birch and 
Swinnerton-Dyer conjecture for abelian varieties over a global field 
of positive characteristic, in \cite{bauer}, \S 4 and \cite{KT}, \S 2 for instance. Our initial interest 
in Theorem~\ref{thm AV} was to try to extend it to complexes of tori 
in the function field case, following the same method as in the number 
field case \cite{demimrn}. Such a generalization should then provide 
results (known in the number field case) about weak and strong approximation 
for linear algebraic groups defined over a global field of positive 
characteristic.

\smallskip

However, as K.~\v{C}esnavi\v{c}ius pointed out to us\footnote{In particular, 
he observed that the analogue of \cite{MilADT}, Prop.~III.0.4.c) is by no
means obvious when henselizations are replaced by completions. 
This analogue is 
actually false without additional assumptions, as shown by T.~Suzuki in 
\cite{suzuki}, Rem~2.7.9}, 
it seems 
necessary to add details to the proof in \cite{MilADT}, 
sections III.3.~and III.8, for two reasons: 

\begin{itemize}

        \item the functoriality of flat cohomology with compact support and the commutativity of several diagrams are not explained in \cite{MilADT}. Even 
in the case of an imaginary number field, a definition of 
$H^r_c(U,{\mathcal F})$ as $H^r(X, j_! {\mathcal F})$ for a fppf 
sheaf ${\mathcal F}$ (which works for the 
\'etale Artin--Verdier Theorem) would not be the right one, 
because it does not provide the key exact sequence \cite{MilADT} 
Prop.~III.0.4.a) in the fppf setting: indeed the proof of this exact sequence 
relies on \cite{MilADT}, Lemma II.2.4., which in turn uses \cite{MilADT}, 
Prop~II.1.1; but the analogue of the latter does not stand anymore with 
\'etale cohomology replaced by fppf cohomology, see also 
Remark \ref{rem diff compl} of the present paper.

It is therefore necessary to work 
with an adhoc definition of compact support cohomology as 
in loc.~cit., \S III.0. Since this definition involves mapping cones, 
commutativities of some 
diagrams have to be checked in the category of complexes and not in 
the derived category (where there is no good functoriality for the 
mapping cones). Typically, the isomorphisms that compute 
$C^{\bullet}(b)$, $C^{\bullet}(b \circ a)$ and 
$C^{\bullet}(c \circ b \circ a)$ in loc.~cit., Prop.~III.0.4.c) 
are not canonical a priori. Hence the required compatibilities 
in loc.~cit., proof of Theorem~III.3.1.~and Lemma~III.8.4.~have to be 
checked carefully.

        \item in the positive characteristic case, it is necessary
(as explained in \cite{MilADT}, \S III.8.) to work with a definition 
of cohomology with compact support involving 
{\it completions} of the local rings of points in $X \setminus U$
instead of their henselizations. The reason is that a local duality 
statement (loc.~cit., Th.~III.6.10) is needed and this one  
only works in the context of complete valuation fields, 
in particular because the $H^1$ groups 
involved have to be locally compact (so that Pontryagin duality makes sense).
It turns out that some 
properties of compact support cohomology (in particular loc.~cit., 
Prop.~III.0.4.c)) are more difficult to establish in this context: 
for example the comparison between cohomology of 
the completion $\widehat{\calo_v}$ and of the henselization $\calo_v$ is not 
as straightforward as in the \'etale case.

\end{itemize}

The goal of this article is to present a detailed proof of 
Theorem~\ref{thm AV} with special regards to the two issues listed 
above. Section~\ref{sect1} is devoted to general properties of fppf 
cohomology with compact support (Prop.~\ref{prop supp compact}), 
which involves some 
homological algebra (Lemma~\ref{lem complex}) 
as well as comparison statements between 
cohomology of $\calo_v$ and $\widehat{\calo_v}$ (Lemma~\ref{lem artin}); 
besides, 
we make the link to classical \'etale cohomology with compact 
support (Lemma~\ref{lem comparison etale}).

We also define a natural topology on the fppf compact support cohomology
groups (see section \ref{section topo})
and prove its basic properties. In section~\ref{sect2}, 
we follow the method of \cite{MilADT}, \S III.8.~to prove Theorem~\ref{thm AV} 
in the function field case. As a corollary, we get a finiteness statement
 (Cor.~\ref{corh2}), which apparently has not been observed before this 
paper. 
The case of a number field 
is simpler once the functorial properties of section~\ref{sect1} have been 
proved; it is treated in section~\ref{sect3}.
Finally, we include two useful
results in homological algebra in an appendix (section~\ref{appsect}).

\smallskip

One week after the first draft of this article was released, Takashi Suzuki 
kindly informed us that in his preprint \cite{suzuki}, he 
obtained (essentially at the same time as us) 
fppf duality results similar to Theorem~\ref{thm AV} in a slightly 
more general context.\footnote{Note, however, that there is still some work 
to do to obtain our Theorem~\ref{thm AV} from the very general Th~3.1.3.~of \cite{suzuki}; compare with section~4.2.~of loc.~cit., where a similar 
task is fullfilled for abelian schemes instead of finite group schemes.} 
His methods are somehow more involved than ours, they use the {\it rational 
\'etale site}, which he developed in earlier papers.

\bigskip

{\bf Notation.} Let $X$ be either a smooth projective curve over a finite 
field $k$ of characteristic $p$,
or the spectrum of the ring of integers $\calo_K$ of a number field $K$. 
Let $K := k(X)$ be the function field of $X$.
Throughout the paper, schemes $S$ are endowed with a big fppf site $(\textup{Sch}/S)_{{\rm fppf}}$ in the sense of \cite[\href{http://stacks.math.columbia.edu/tag/021R}{Tag 021R}, \href{http://stacks.math.columbia.edu/tag/021S}{Tag 021S}, \href{http://stacks.math.columbia.edu/tag/03XB}{Tag 03XB}]{SP}. By construction, the underlying category in $({\rm Sch}/S)_{\rm fppf}$ is small and the family of coverings for this site is a set. The corresponding topos is independent 
of the choices made thanks 
to \cite[\href{http://stacks.math.columbia.edu/tag/00VY}{Tag 00VY}]{SP}.
In contrast with \cite{SGA4}, the construction of the site $({\rm Sch}/S)_{\rm fppf}$ in \cite{SP} does not require the existence of universes. The reader who is ready to accept this axiom can replace the site $({\rm Sch}/S)_{\rm fppf}$ by the big fppf site from \cite{SGA4}.

Unless stated otherwise, cohomology is fppf cohomology with respect to this site.

For any closed point $v \in X$, let $\mathcal{O}_v$ (resp.~$\widehat{\mathcal{O}_v}$) be the henselization (resp.~the completion) of the local ring 
$\calo_{X, \, v}$ 
of $X$ at $v$. Let $K_v$ (resp.~$\widehat{K_v}$) be the fraction field of $\mathcal{O}_v$  (resp.~$\widehat{\mathcal{O}_v}$). Let $U$ be a non empty 
Zariski open subset of $X$ and denote by $j : U  \to X$ the 
corresponding open immersion.
By \cite{matsu}, \S 34, the local ring $\calo_{X, \, v}$ of $X$ at $v$ 
is excellent (indeed $\calo_{X, \, v}$ is either of mixed characteristic
or the localization of a ring of finite type over a field); hence 
so are $\mathcal{O}_v$ (by \cite{ega4}, Cor.~18.7.6) 
as the  henselization of an excellent ring, 
and $\widehat{\mathcal{O}_v}$ as a complete Noetherian local ring 
(\cite{matsu}, \S 34).

The piece of notation ``$v \not \in U$'' means that
we consider all places $v$ corresponding to closed points of $X \setminus U$
{\it plus the real places in the number field case}. If $v$ is a real
place, we set $K_v=\widehat{K_v}=\mathcal{O}_v=\widehat{\mathcal{O}_v}$
for the completion of $K$ at $v$, and we denote 
by $H^*(K_v,M)$ the \emph{Tate (or modified)} cohomology groups of a $\textup{Gal}(\overline{K_v}/K_v)$-module $M$.

If $\mathcal{F}$ is a fppf sheaf of abelian groups 
on $U$, define the Cartier dual 
$\mathcal{F}^D$ to be the fppf sheaf 
$\mathcal{F}^D := \underline{\Hom}(\mathcal{F}, \G)$. Notation 
as $\Gamma(U,\mathcal{F})$ stands for the group of sections 
of $\mathcal{F}$ over $U$, and $\Gamma_Z(U,\mathcal{F})$ for the group 
of sections that vanish over $U \setminus Z$. 
If $E$ is a field (e.g.~$E=K_v$ or $E=\widehat{K_v}$) and 
$i : \spec E \to U$ is an $E$-point of $U$, we will frequently 
write $H^r(E,\mathcal{F})$ for $H^r(\spec E,i^* \mathcal{F})$. Similarly 
for an open subset $V \subset U$, the piece of notation 
$H^r(V,\mathcal{F})$ (resp.~$H^r_c(V,\mathcal{F})$) stands for 
$H^r(V,\mathcal{F}_{\vert V})$ (resp.~$H^r_c(V,\mathcal{F}_{\vert V})$).

A finite group scheme $N$ over a field $E$ of characteristic 
$p>0$ is {\it local} (or equivalently {\it infinitesimal}, as in \cite{demgab}, II.4.7.1) if it is connected (in particular this implies 
$H^0(E',N)=0$ for every field extension $E'$ of $E$). 
Examples of such group schemes
are $\mu_p$ (defined by the affine equation $y^p=1$) and $\alpha_p$ 
(defined by the equation $y^p=0$).  

Let $S$ be an $\F_p$-scheme.
A finite $S$-group scheme 
$N$ is of {\it height 1} if the relative Frobenius map 
$F_{N/S}$ (cf.~\cite{MilADT}, \S III.0) is trivial. 

For any topological abelian group $A$, let $A^* := \Hom_{\textup{cont.}}(A, \Q/\Z)$ be the group of continuous homomorphisms from $A$ to $\Q/\Z$ 
(where $\Q/\Z$ is considered as a discrete group) equipped with 
the compact--open topology. 
A morphism $f : A \to B$ 
of topological groups is {\it strict} if it is continuous, and 
the restriction $f : A \to f(A)$ is 
an open map (where the topology on $f(A)$ is induced by $B$). This is 
equivalent to saying that $f$ induces an isomorphism of the topological 
quotient $A/\ker f$ with the topological subspace $f(A) \subset B$.

Concerning sign conventions in homological algebra, we tried to follow the conventions in \cite{SP} throughout the text.

\section{Fppf cohomology with compact support} \label{sect1}

Define $Z := X \setminus U$ and $Z' := \coprod_{v \in Z} \spec(\widehat{K_v})$ (disjoint union).
Then we have a natural morphism $i : Z' \to U$.
Let $\mathcal{F}$ be a sheaf of abelian groups 
on $({\rm Sch} / U)_{\rm fppf}$. Let $I^\bullet(\mathcal{F})$ be an injective resolution of $\mathcal{F}$ over $U$. 
Denote by $\mathcal{F}_v$ and $I^\bullet(\mathcal{F})_v$ their respective 
pullbacks to $\spec K_v$, for $v \notin U$.

Given a morphism of schemes $f : T \to S$, the fppf pullback functor $f^*$ is exact (see \cite[\href{http://stacks.math.columbia.edu/tag/021W}{Tag 021W}, \href{http://stacks.math.columbia.edu/tag/00XL}{Tag 00XL}]{SP}) and it admits an exact left adjoint $f_!$ (see \cite[\href{http://stacks.math.columbia.edu/tag/04CC}{Tag 04CC}]{SP}), hence $f^*$ maps injective (resp.~flasque) 
objects to injective (resp.~flasque) objects. 
Therefore $I^\bullet(\mathcal{F})_v$ is an injective resolution 
of $\mathcal{F}_v$.

\smallskip

As noticed by A.~Schmidt, the definition of the modified fppf 
cohomology groups in the number field case in \cite{MilADT}, III.0.6 (a), 
has to be written more precisely, because of the non canonicity of the 
mapping cone in the derived category. We are grateful to him for the 
following alternative definition. 

Let $\Omega_\R$ denote the set of 
real places of $K$. For $v \in \Omega_\R$, 
let 
$a^v : ({\rm Sch}/\spec({K}_v))_{\textup{fppf}} 
\to \spec({K}_v)_{\textup{\'et}}$ be the natural morphism of sites, where $S_{\textup{\'et}}$ denotes the small \'etale site on a scheme $S$. Since $K_v$ is a perfect field, the direct image functor $a^v_*$  associated to $a^v$ is exact. Hence, by \cite{SGA4}, V, Remark~4.6 and Prop.~4.9, 
the functor $a^v_*$ maps $I^\bullet(\mathcal{F})_v$ to a flasque resolution 
$a^v_* I^\bullet(\mathcal{F})_v$ of $a^v_* \mathcal{F}_v$. Following \cite{GS} \S 2, there is a natural acyclic resolution 
$D^\bullet(a^v_* \mathcal{F}_v) \to a^v_* \mathcal{F}_v$ 
of the $\textup{Gal}(\overline{{K_v}}/{K_v}) = \Z / 2\Z$-module 
$a^v_* \mathcal{F}_v$ 
(identified with $\mathcal{F}_v(\spec(\overline{{K_v}}))$). 
Splicing the resolutions $D^\bullet(a^v_* \mathcal{F}_v)$ 
and $a^v_* I^\bullet(\mathcal{F})_v$ together, 
one gets a complete acyclic resolution $\widehat{I}^\bullet(\mathcal{F}_v)$ 
of the $\textup{Gal}(\overline{{K_v}}/{K_v})$-module 
$a^v_* \mathcal{F}_v$, which computes the Tate cohomology of  
$a^v_* \mathcal{F}_v$. And by construction, 
there is a natural morphism 
$\widehat{i}_v : a^v_* I^\bullet(\mathcal{F})_v 
\to \widehat{I}^\bullet(\mathcal{F}_v)$.

\smallskip

As suggested by \cite{MilADT}, section III.0, define $\Gamma_c(U, I^\bullet(\mathcal{F}))$ to be the following object in the category of complexes of abelian groups:
$$\Gamma_c(U, I^\bullet(\mathcal{F})) := \Cone\left(\Gamma(U, I^\bullet(\mathcal{F})) \to \Gamma(Z', i^* I^\bullet(\mathcal{F})) \oplus \bigoplus_{v \in \Omega_\R} \Gamma({K_v}, \widehat{I}^\bullet(\mathcal{F}_v))  \right)[-1] \, ,$$
and $H^r_c(U, \mathcal{F}) := H^r(\Gamma_c(U, I^\bullet(\mathcal{F})))$. We will also denote by $R\Gamma_c(U,\mathcal{F})$ the complex $\Gamma_c(U, I^\bullet(\mathcal{F}))$ viewed in the derived category of fppf sheaves.
Observe that in the number field case the groups $H^r_c(U, \mathcal{F})$
may be non zero even for negative $r$.
In the function field case we have $H^r_c(U, \mathcal{F})=0$ for $r <0$, 
and also (by Proposition~\ref{prop supp compact}
below) $H^0_c(U, \mathcal{F})=0$ if we assume further $U \neq X$ 
(the map $H^0(U, \mathcal F) \to H^0({\widehat K}_v,\mathcal F)$
being injective for each $v \not \in U$).

\smallskip

From now on, we will abbreviate $\Cone(...)$ by $C(...)$. 

\begin{prop} \label{prop supp compact}

\hfill \begin{enumerate}
	\item \label{prop supp compact-1} Let $\mathcal{F}$ be a 
sheaf of abelian groups on $U_{{\rm fppf}}$.
There is a natural exact sequence, for all $r \geq 0$, 
$$\dots \to H^r_c(U, \mathcal{F}) \to H^r(U, \mathcal{F}) \to \bigoplus_{v \not \in U} H^r(\widehat{K}_v, \mathcal{F}) \to H^{r+1}_c(U, \mathcal{F}) \to \dots \, .$$
	\item \label{prop supp compact-2} For any short exact sequence 
$$0 \to \mathcal{F'} \to \mathcal{F} \to \mathcal{F}'' \to 0$$
of sheaves of abelian groups on $U$, there is a long exact sequence
$$\dots \to H^r_c(U, \mathcal{F}') \to H^r_c(U, \mathcal{F}) \to H^r_c(U, \mathcal{F}'') \to H^{r+1}_c(U, \mathcal{F}') \to \dots \, .$$

	\item \label{prop supp compact-3} For any flat affine commutative group scheme $\mathcal{F}$ of finite type over $U$, and any non empty 
open subscheme $V \subset U$, there is a canonical exact sequence
$$\dots \to H^r_c(V, \mathcal{F}) \to H^r_c(U, \mathcal{F}) \to \bigoplus_{v \in U \setminus V} H^r(\widehat{\mathcal{O}}_v, \mathcal{F}) \to H^{r+1}_c(V, \mathcal{F}) \to \dots \, ,$$
and the following natural diagram commutes:
\begin{changemargin}{-2cm}{2cm}
\[
\xymatrix{
& \bigoplus_{v \notin V} H^{r-1}(\widehat{K_v},\mathcal{F}) \ar[d] & \bigoplus_{v \notin U} H^{r-1}(\widehat{K_v}, \mathcal{F}) \ar[l]_-{i_2} \ar[d] & \\
\bigoplus_{v \in U \setminus V} H^{r-1}(\widehat{\mathcal{O}_v}, \mathcal{F})
 \ar[r] \ar[ru]^-{i_1} & H^r_c(V, \mathcal{F}) \ar[r] \ar[d] & H^r_c(U, \mathcal{F}) \ar[r] \ar[d] & \bigoplus_{v \in U \setminus V} H^{r}(\widehat{\mathcal{O}_v}, \mathcal{F}) \\
& H^r(V,\mathcal{F}) \ar[d]& H^r(U,\mathcal{F}) \ar[d] \ar[l]_-{\textup{Res}} \ar[ru] & \\
& \bigoplus_{v \notin V} H^r(\widehat{K_v},\mathcal{F}) \ar[r]^-{\pi} & \bigoplus_{v \notin U} H^r(\widehat{K_v}, \mathcal{F}) \, , &
}
\]
\end{changemargin}
where $i_1$ (resp.~$i_2$) is obtained by puting  
$0$ at the places $v \not 
\in U$ (resp.~$v \in U\setminus V$) and $\pi$ is the natural projection. 

	\item \label{prop supp compact-4} If $\mathcal{F}$ is represented 
by a smooth group scheme, then $H^r_c(U, \mathcal{F}) \cong H^r_{\et,c}(U, \mathcal{F})$ for $r \neq 1$, where $H_{\et,c} ^*$ stands for modified \'etale cohomology with compact support as defined in \cite{GS}, \S 2.
In particular for such $\mathcal{F}$ we have 
$H^r_c(U,\mathcal{F}) \cong H^r_{\et}(X, j_! \mathcal{F})$ in the 
function field case.
 If in addition the generic fiber $\mathcal{F}_K$ is a finite $K$-group scheme, then $H^1_c(U, \mathcal{F}) \cong H^1_{\et,c}(U, \mathcal{F})$ 
(which is identified with 
$H^1_{\et}(X, j_! \mathcal{F})$ in the function field case).

\end{enumerate}
\end{prop}

\begin{rema} \label{rem diff compl}
{\rm Unlike what happens in \'etale cohomology, 
the groups $H^1(\calo_v,{\cal F})$ and 
$H^1(\widehat{\calo_v},{\cal F})$ cannot in general be identified with 
the group $H^1(k(v),F(v))$, where $k(v)$ is the residue 
field at $v$ and $F(v)$ the fiber of ${\cal F}$ over $k(v)$.
For example this already fails for ${\cal F}=\mu_p$ and 
$\widehat{\calo_v}=\F_p[[t]]$, because 
by the Kummer exact sequence
$$0 \to \mu_p \to \G \stackrel{.p}{\to} \G \to 0 $$
in fppf cohomology, the group 
$H^1(\widehat{\calo_v},{\cal F})=\widehat{\calo_v}^*/\widehat{\calo_v}^{*^p}$ is an infinite dimensional $\F_p$-vector space, while $H^1(k(v),F(v))=k(v)^*/k(v)^{*^p}=0$.
The situation is better for $r \geq 2$ by \cite{Toe}, 
Cor.~3.4: namely the natural maps from 
$H^r(\calo_v,{\cal F})$ and $H^r(\widehat{\calo_v},{\cal F})$ to 
$H^r(k(v),F(v))$ are isomorphisms.
}
\end{rema}

\smallskip

Before proving Proposition \ref{prop supp compact}, 
we need the following lemmas.
We start with a lemma in homological algebra:
\begin{lem} \label{lem complex}
Let $\mathcal{A}$ be an abelian category with enough injectives and let ${\bf C}(\mathcal{A})$ (resp.~${\bf D}(\mathcal{A})$) denote the category (resp.~the derived category) of bounded below cochain complexes in $\mathcal{A}$. Consider a commutative diagram in ${\bf C}(\mathcal{A})$:
\[
\xymatrix{
A \ar[d]_f \ar[r]^-{\alpha} & B \oplus E \ar[d]^-{(\id, g)} \\
A' \ar[r]^-{\alpha'} & B \oplus E' \, ,
}
\]
and denote by $\pi_B$ (resp.~$\pi'_B$) the projection $B \oplus E \to B$ (resp.~$B \oplus E' \to B$).

Assume that the natural morphism $C(f) \to C(g)$ in ${\bf C}(\mathcal{A})$ is a quasi-isomorphism. Then there exists a canonical commutative diagram in ${\bf D}(\mathcal{A})$:
\[
\xymatrix{
(B \oplus E')[-1] \ar[d] & B[-1] \ar[l]_-{i'_{B}} \ar[d] & B \oplus E \ar[r]^-{(\id, g)} & B \oplus E' \ar[d] \\
C(\alpha')[-1] \ar[r] \ar[d] & C(\pi_B \circ \alpha)[-1] \ar[r] \ar[d] & E \ar[r] \ar[u]^-{i_E} & C(\alpha') \\
A' \ar[d]^-{\alpha'} & A \ar[d]^-{\pi_B \circ \alpha} \ar[l]_{f} \ar[r]^-{\alpha} & B \oplus E \ar[u]^-{\pi_E} & \\
B \oplus E' \ar[r]^-{\pi'_{B}} & B \, , & & 
}
\]
where the second row and the first two columns are exact triangles.
\end{lem}

\dem{}
The assumption that $C(f) \to C(g) \cong C({\rm Id} \oplus g)$ 
is a quasi-isomorphism implies that $C(\alpha) \to C(\alpha')$ is a quasi-isomorphism (see for instance Proposition 1.1.11 in \cite{BBD} or Corollary A.14 in \cite{MHS}).

Functoriality of the mapping cone in the category ${\bf C}(\mathcal{A})$ 
gives the following 
diagram in ${\bf C}(\mathcal{A})$, where the second row (by \cite{MilADT}, 
Prop.~II.0.10, or \cite{KS}, proof of Theorem 11.2.6) and the columns are exact triangles in the derived category:
\[
\xymatrix{
& (B\oplus E)[-1] \ar[dl]_-{(\id, g)} \ar[d] & B[-1] \ar[l]_-{i_B} \ar[d] \ar[r]^= & B[-1] \ar[d] & & \\
(B \oplus E')[-1] \ar[d] & C(\alpha)[-1] \ar[ld] \ar[d] \ar[r] & C(\pi_B \circ \alpha)[-1] \ar[d] \ar[r] & C(\pi_B)[-1] \ar[d] \ar[r] \ar[dr] \ar@{}|(.5){\star}[drr] & C(\alpha) & \\
C(\alpha')[-1] \ar[d] & A \ar[ld]_-f \ar[d]^{\alpha} \ar[r]^-= & A \ar[d]^-{\pi_B \circ \alpha} \ar[r]^-\alpha & B \oplus E \ar[d]^-{\pi_B} \ar[r]^-{\pi_E} & E \ar[r]^-{i_E} & B \oplus E \ar[ul] \\
A' \ar[d]^-{\alpha'} & B \oplus E \ar[dl]_-{(\id,g)} \ar[r]^-{\pi_B} & B \ar[r]^-= & B & & \\
B \oplus E' & & & & \, .
}
\]
As usual, notation as $\pi_B$, $\pi_E$ denotes projections and $i_B$, $i_E$ 
are given by puting $0$ at the missing piece. Note also that due to our sign conventions, the horizontal map $C(\pi_B)[-1] \to C(\alpha)$ is given by the natural map with a $(-1)$-sign. 

This diagram is commutative in ${\bf C}(\mathcal{A})$, except the 
square $\star$ 
which is commutative up to homotopy. 
Indeed, this square defines two maps $f, g : C(\pi_B)[-1] \to C(\alpha)$, which are given in degree $n$ by two maps $f^n, g^n : B^{n-1} \oplus (B^n \oplus E^n) \to (B^n \oplus E^n) \oplus A^{n+1}$, where $f^n(b',b,e) := -(b,e,0)$ and $g^n(b',b,e) := -(0,e,0)$. Consider now the maps $s^n : B^{n-1} \oplus (B^n \oplus E^n) \to (B^{n-1} \oplus E^{n-1}) \oplus A^{n}$ defined by
$s^n(b',b,e) := (b',0,0)$. 
Then the collection $(s^n)$ is a homotopy between $f$ and $g$. Hence the square $\star$ is commutative up to the homotopy $(s^n)$.

Since the map $C(\alpha) \to C(\alpha')$ is a quasi-isomorphism, and since the natural map $C(\pi_B)[-1] \to E$ is a homotopy equivalence, the lemma follows from the commutativity and the exactness of the previous diagram.
\enddem

We now need the following result, for which we did not 
find a suitable reference:
\begin{lem} \label{lem dense H1}
 Let $A$ be a henselian valuation ring with fraction field $K$.
 Let $\widehat A$ be the completion of $A$ for the valuation topology and 
$\widehat K:=\fraco \widehat A$.
Assume that $\widehat{K}$ is separable over $K$.
\begin{enumerate}
	\item Let $G$ be a $K$-group scheme locally of finite type.
Then the map $H^1(K, G) \to H^1(\widehat{K}, G)$ has dense image.

	\item Assume that $\widehat{A}$ is henselian. 
Let $G$ be a flat $A$-group scheme locally of finite presentation.
Then the map $H^1(A, G) \to H^1(\widehat{A}, G)$ has dense image.
\end{enumerate}
\end{lem}

Here the topology on the pointed sets $H^1(\widehat{A}, G)$ and 
$H^1(\widehat{K}, G)$ is provided by \cite{Ces},~\S~3.

\begin{rema}
{\rm \hfill
\begin{itemize}
	\item The assumption that $\widehat{K}$ is separable over $K$ is 
satisfied if $A$ is an excellent discrete valuation ring. 
	\item In the second statement, the assumption that $\widehat{A}$ 
is henselian is satisfied if the valuation on $A$ has height $1$ (special 
case of \cite{riben}, section F, Th.~4). 
This assumption is used in the proof below to apply \cite{Ces}, Theorem B.5.
Note also that in general, $\widehat{A}$ is not the same as the completion of $A$ for the 
$\mathfrak{m}$-adic topology (where $\mathfrak{m}$ denotes the maximal ideal of 
$A$).
\end{itemize}
}
\end{rema}

\dem{ of Lemma~\ref{lem dense H1}}
We prove both statements at the same time. 
Let $E$ be either $A$ or $K$, set $S=\spec E$.
 Let ${\bf B}G$ denote the classifying $E$-stack of $G$-torsors. We need to prove that ${\bf B}G(E)$ is dense in ${\bf B}G(\widehat{E})$.
It is a classical fact that ${\bf B}G$ is an algebraic stack
(\cite[\href{https://stacks.math.columbia.edu/tag/0CQJ}{Tag 0CQJ}]{SP} 
and \cite[\href{https://stacks.math.columbia.edu/tag/06PL}{Tag 06PL}]{SP}).
 Let $x \in {\bf B}G (\widehat{E})$ and $U \subset {\bf B}G (\widehat{E})$ 
be an open subcategory (in the sense of \cite{Ces}, 2.4) containing $x$. 
We need to find an object $x' \in {\bf B}G (E)$ that maps to 
$U \subset {\bf B}G (\widehat{E})$.
 Using \cite{Ces}, Theorem B.5 and Remark B.6 (applied to the $S$-scheme 
$\spec R:=\spec \widehat{E}$), there exists an affine 
scheme $Y$, a smooth $S$-morphism $\pi : Y \to {\bf B} G$ and $y \in 
Y(\widehat{E})$ such that $\pi_{\widehat{E}}(y) = x$, where $\pi_{\widehat{E}} : Y(\widehat E) \to BG(\widehat E)$ is the map induced by $\pi$.
In particular, $Y \to S$ is smooth because so are $\pi$ and 
${\bf B}G \to S$ (the latter by \cite{Ces}, Prop.~A.3). Hence
$Y$ is locally of finite presentation over $S$.
By assumption, $\pi_{\widehat{E}}^{-1}(U) \subset Y(\widehat{E})$ is 
an open subset containing $y$. Hence \cite{MB}, Corollary 1.2.1 (in the 
discrete valuation ring case, 
it is Greenberg's approximation Theorem) implies that 
$Y(E) \cap \pi_{\widehat{E}}^{-1}(U) \neq \emptyset$. Applying $\pi_E$, 
we get that the image of ${\bf B}G (E)$ meets $U$,
which proves the required result.
\enddem

The previous lemma is useful to prove the following crucial (in the 
function field case) statement. 
For a local integral domain $A$ with maximal
ideal $\mathfrak{m}$, fraction field $K$ and residue field $\kappa$, 
and $\mathcal{F}$ an fppf sheaf on $\spec A$ with an injective 
resolution $I^\bullet(\mathcal{F})$, define 
$$\Gamma_{\mathfrak{m}}(A, I^\bullet(\mathcal{F})) := 
\Cone\left(\Gamma(\spec A, I^\bullet(\mathcal{F})) \rightarrow 
\Gamma(\spec K,I^\bullet(\mathcal{F})) \right)[-1]$$ and 
$H^r_{\mathfrak{m}}(A, \mathcal{F}) := 
H^r(\Gamma_{\mathfrak{m}}(A, I^\bullet(\mathcal{F})))$ (the 
cohomology with support in $\spec \kappa$). 
We have a localization long exact sequence 
(\cite{MilADT}, Prop.~III.0.3)
$$...\to H^r_{\mathfrak{m}}(A, \mathcal{F}) \to H^r(A, \mathcal{F}) \to 
H^r(K,\mathcal{F}) \to H^{r+1}_{\mathfrak{m}}(A, \mathcal{F}) \to ...$$

\begin{lem} \label{lem artin}
Let $A$ be an excellent henselian discrete valuation ring, with maximal ideal $\mathfrak{m}$. Let $\mathcal{F}$ be a flat affine commutative group scheme 
of finite type over $\spec A $.
Then for all $r \geq 0$, the morphism $H^r_{\mathfrak{m}}(A, \mathcal{F}) \to H^r_{\mathfrak{m}}(\widehat{A}, \mathcal{F})$ is an isomorphism.
\end{lem}

\begin{rema} \label{rem acyclic}
{\rm Let $I^\bullet(\mathcal{F})$ be an injective resolution of 
$\mathcal{F}$ viewed as an fppf sheaf.
Another formulation of Lemma~\ref{lem artin} is that the natural morphism 
$\Gamma_{\mathfrak{m}}(A,I^\bullet(\mathcal{F})) 
\to \Gamma_{\mathfrak{m}}(\widehat{A}, I^\bullet(\mathcal{F}))$ is an 
isomorphism in the derived category.
The injective resolution $I^\bullet(\mathcal{F})$ can be replaced by any 
complex of flasque fppf sheaves that is quasi-isomorphic to $\mathcal{F}$
(indeed the fppf pullback functor $f^*$ associated to 
$f : \spec \widehat A \to \spec A$ sends flasque resolutions to 
flasque resolutions, because $f^*$ is exact and preserves flasque 
sheaves). 

\smallskip

Also note that Lemma~\ref{lem artin} 
is a variant of \cite{suzuki}, Prop.~2.6.2: our result is slightly more general in the affine case, while the notion of cohomological approximation in \cite{suzuki} is a priori a little stronger than the conclusion of Lemma~\ref{lem artin}. In addition, this lemma answers a variant of a question raised after
Prop~2.6.2.~of loc.~cit.~(under a flatness assumption). }
\end{rema}

\dem{ of Lemma~\ref{lem artin}}

\begin{itemize}
	\item $r = 0$: 

Since $\mathcal{F}$ is separated (as an affine scheme), the morphisms $H^0(A, \mathcal{F}) \to H^0(K, \mathcal{F})$ and $H^0(\widehat{A}, \mathcal{F}) \to H^0(\widehat{K}, \mathcal{F})$ are injective, which implies that
$$H^0_{\mathfrak{m}}(A, \mathcal{F}) = H^0_{\mathfrak{m}}(\widehat{A}, \mathcal{F}) = 0 \, .$$
	
	\item $r=1$:

Consider the following commutative diagram with exact rows:
\begin{equation} \label{diag artin 1}
\begin{gathered}
\resizebox{12,5cm}{!}{
\xymatrix{
H^0(A, \mathcal{F}) \ar[r] \ar[d] & H^0(K, \mathcal{F}) \ar[r] \ar[d] & H^1_{\mathfrak{m}}(A, \mathcal{F}) \ar[r] \ar[d] & H^1(A, \mathcal{F}) \ar[r] \ar[d] & H^1(K, \mathcal{F}) \ar[d] \\
H^0(\widehat{A}, \mathcal{F}) \ar[r] & H^0(\widehat{K}, \mathcal{F}) \ar[r] & H^1_{\mathfrak{m}}(\widehat{A}, \mathcal{F}) \ar[r] & H^1(\widehat{A}, \mathcal{F}) \ar[r] &  H^1(\widehat{K}, \mathcal{F})\, .}}
\end{gathered}
\end{equation}

Since $A$ is excellent, Artin approximation (see 
\cite{Art}, Theorem 1.12) implies that the morphism 
$H^1(A, \mathcal{F}) \to H^1(\widehat{A}, \mathcal{F})$ is injective: indeed, given a ($\spec A$)-torsor $\mathcal{P}$ under $\mathcal{F}$, $\mathcal{P}$ is locally of finite presentation, and Artin approximation ensures that $\mathcal{P}(\widehat{A}) \neq \emptyset$ implies that $\mathcal{P}(A) \neq \emptyset$.

The affine $A$-scheme of finite type ${\mathcal F}$ is of the form 
$\spec(A[x_1,...,x_n]/(f_1,...,f_r))$, where $f_1,...,f_r$ are
polynomials. Since the discrete valuation ring $A$ satisfies 
$A=K \cap \widehat A \subset \widehat K$, 
the square on the left hand side in \eqref{diag artin 1} is cartesian.

Hence an easy diagram chase implies that $H^1_{\mathfrak{m}}(A, \mathcal{F})  \to H^1_{\mathfrak{m}}(\widehat{A}, \mathcal{F})$ is injective.

By Proposition A.6 in \cite{GP}, the right hand side square in \eqref{diag artin 1} is cartesian.
In addition, $H^0(\widehat{A}, \mathcal{F}) \subset H^0(\widehat{K}, \mathcal{F})$ is open (\cite{GGMB}, Prop.~3.3.4), and $H^0(K, \mathcal{F}) \subset H^0(\widehat{K}, \mathcal{F})$ is dense by \cite{GGMB}, Proposition 3.5.2 (weak approximation for $\mathcal{F}$).

Therefore, an easy diagram chase implies that the map $H^1_{\mathfrak{m}}(A, \mathcal{F})  \to H^1_{\mathfrak{m}}(\widehat{A}, \mathcal{F})$ is surjective.

	\item $r=2$:

Consider the commutative diagram with exact rows:
\begin{equation} \label{diag artin 2}
\begin{gathered}
\resizebox{12,5cm}{!}{
\xymatrix{
H^1(A, \mathcal{F}) \ar[r] \ar[d] & H^1(K, \mathcal{F}) \ar[r] \ar[d] & H^2_{\mathfrak{m}}(A, \mathcal{F}) \ar[r] \ar[d] & H^2(A, \mathcal{F}) \ar[r] \ar[d] & H^2(K, \mathcal{F}) \ar[d] \\
H^1(\widehat{A}, \mathcal{F}) \ar[r] & H^1(\widehat{K}, \mathcal{F}) \ar[r] & H^2_{\mathfrak{m}}(\widehat{A}, \mathcal{F}) \ar[r] & H^2(\widehat{A}, \mathcal{F}) \ar[r] & H^2(\widehat{K}, \mathcal{F}) \, .
}}
\end{gathered}
\end{equation}

By \cite{Toe}, Corollary 3.4, the map $H^2(A, \mathcal{F}) \to H^2(\widehat{A}, \mathcal{F})$ is an isomorphism. And we already explained (in the case $r=1$) that the left hand side square in \eqref{diag artin 2} is cartesian. Hence a diagram chase proves that the map $H^2_{\mathfrak{m}}(A, \mathcal{F}) \to H^2_{\mathfrak{m}}(\widehat{A}, \mathcal{F})$ is injective.

Using \cite{GGMB}, Proposition 3.5.3.(3), the map $H^2(K, \mathcal{F}) \to H^2(\widehat{K}, \mathcal{F})$ is also an isomorphism.
By \cite{Ces}, Proposition 2.9 (e), the map $H^1(\widehat{A}, \mathcal{F}) \to H^1(\widehat{K}, \mathcal{F})$ is open.
Lemma \ref{lem dense H1} implies that the map $H^1(K, \mathcal{F}) \to H^1(\widehat{K}, \mathcal{F})$ has dense image. By diagram chase, we get that the map $H^2_{\mathfrak{m}}(A, \mathcal{F}) \to H^2_{\mathfrak{m}}(\widehat{A}, \mathcal{F})$ is surjective.

	\item $r \geq 3$:

Corollary 3.4 in \cite{Toe} implies that the morphisms $H^{r-1}(A, \mathcal{F}) \to H^{r-1}(\widehat{A}, \mathcal{F})$ and $H^r(A, \mathcal{F}) \to H^r(\widehat{A}, \mathcal{F})$ are isomorphisms. Proposition 3.5.3.(3) in \cite{GGMB} implies that the maps $H^{r-1}(K, \mathcal{F}) \to H^{r-1}(\widehat{K}, \mathcal{F})$ and $H^r(K, \mathcal{F}) \to H^r(\widehat{K}, \mathcal{F})$ are isomorphisms. Therefore, the five-lemma proves that $H^r_{\mathfrak{m}}(A, \mathcal{F}) \to H^r_{\mathfrak{m}}(\widehat{A}, \mathcal{F})$ is an isomorphism. \enddem
\end{itemize}

\begin{rema}
{\rm We will apply the previous lemma to a finite and flat commutative
group scheme $N$. As was pointed out to us by K.~\v{C}esnavi\v{c}ius, 
it is then possible to argue without using 
Corollary 3.4 in \cite{Toe} (whose proof is quite involved): indeed 
there exists (cf.~\cite{MilADT}, Th.~III.A.5) an exact sequence 
$$0 \to N \to G_1 \to G_2 \to 0$$
of affine $A$-group schemes such that $G_1$ and $G_2$ are smooth. Now  
for $i >0$ we have $H^i(A,G_j) \cong H^i(\widehat A,G_j)$ ($j=1,2$)
by \cite{MilEC}, Rem.~III.3.11
because $A$ and $\widehat A$ are henselian, and fppf cohomology coincides
with \'etale cohomology for smooth group schemes. It remains to apply 
the five-lemma to get $H^i(A,N) \cong H^i(\widehat A,N)$ for $i \geq 2$, 
which is the input from \cite{Toe} that we used in the proof.
}
\end{rema}

The following lemma is a version of the excision property for fppf cohomology with respect to \'etale morphisms:
\begin{lem} \label{lem excision}
Let $X, X'$ be schemes, $Z \hookrightarrow X$ (resp.~$Z' \hookrightarrow X'$) be closed subschemes, $\pi : X' \to X$ be an \'etale morphism. Assume that $\pi$ restricted to $Z'$ is an isomorphism from $Z'$ to $Z$ and that $\pi(X' \setminus Z') \subset X \setminus Z$. Let $\mathcal{F}$ be a sheaf on $(\textup{Sch}/X)_{\textup{fppf}}$. Then for all $r \geq 0$, the natural morphism $H^r_Z(X, \mathcal{F}) \to H^r_{Z'}(X', \pi^* \mathcal{F})$ is an isomorphism.
\end{lem}

\dem{}
Since $\pi^*$ is exact and maps injective objects to injective objects, 
the proof is exactly the same as the proof
of \cite{MilEC}, Proposition III.1.27.
\enddem

We continue with a lemma comparing the definition of modified 
\'etale cohomology with compact support in \cite{GS} and our 
definition of modified fppf cohomology with compact support. For any scheme $T$, consider the morphisms of sites 
\[
\xymatrix{
({\textup Sch}/T)_{\textup{fppf}} \ar@/_2pc/[rr]^-{a_T} \ar[r]^-{\varepsilon_T} & ({\textup Sch}/T)_{\textup{\'et}} \ar[r]^-{\pi_T} & T_{\et} ,
}
\]
where $({\textup Sch}/T)_{\textup{\'et}}$ denotes the big \'etale site 
of $T$.
Recall that $Z:=X \setminus U$, 
$Z':=\coprod_{v \in Z} \spec ({\widehat K}_v)$, $j : U \to X$ is the open
immersion  and $i : Z' \to U$ is the natural morphism.
Set $a:=a_U$ and $\varepsilon:=\varepsilon_U$.

Let $\mathcal{F}$ be a sheaf on $U_\et$, 
and let $\pi_X^* j_! \mathcal{F} \to J^\bullet(\mathcal{F})$ be an injective 
resolution in the big \'etale topos of $X$. By \cite[\href{http://stacks.math.columbia.edu/tag/0758}{Tag 0758}]{SP} and \cite[\href{http://stacks.math.columbia.edu/tag/04BT}{Tag 04BT}]{SP}, the restriction $J^\bullet(\mathcal{F})_\et$ of $J^\bullet(\mathcal{F})$ to the small \'etale site of $X$ is an injective resolution of $j_! \mathcal{F}$. 
For every place $v \not \in U$ of $K$, let $\mathcal{F}_v$
be the pull-back of $\mathcal{F}$ to $(\spec K_v)_\et$.
As in the fppf case (explained in the beginning of section~\ref{sect1}),
we have for $v$ real 
a complete resolution $\widehat{J}^\bullet(\mathcal{F}_v)$
of the $\gal(\ov K_v/K_v)$-module $\mathcal{F}_v$, which computes its 
Tate cohomology.
Following \cite{GS}, section 2, we define
$$\Gamma_{\et,c}(U,J^\bullet(\mathcal{F})) := \textup{Cone}\left(\Gamma(X, J^\bullet(\mathcal{F})_\et) \to \bigoplus_{v \in \Omega_\R} \Gamma({K_v}, \widehat{J}^\bullet(\mathcal{F}_v))  \right)[-1] \, ,$$
and $H^r_{\et,c}(U, \mathcal{F}) := H^r(\Gamma_{\et,c}(U, J^\bullet(\mathcal{F})))$. 

\smallskip

Denote by $R\Gamma_{\et,c}(U,\mathcal{F})$ the complex 
$\Gamma_{\et,c}(U,J^\bullet(\mathcal{F}))$ (viewed in the derived category 
of abelian groups); similarly for $v$ real, set 
$\widehat{R\Gamma_{\et}}(K_v,\mathcal{F})$ (resp.~$\widehat{R\Gamma}(K_v,a^*\mathcal{F})$) for the 
complex $\Gamma({K_v}, \widehat{J}^\bullet(\mathcal{F}_v))$ 
(resp.~$\Gamma({K_v}, \widehat{I}^\bullet((a^*\mathcal{F})_v))$, 
where $I^{\bullet}(a^*\mathcal{F})$ is a flasque resolution of 
$a^*\mathcal{F}$, cf.~beginning of
section~\ref{sect1}) in the derived
category of \'etale sheaves (resp.~fppf sheaves) over $\spec K_v$. Finally, let $R\Gamma_{\et, Z}(X, j_! \mathcal{F})$ denote the complex 
\[
\Gamma_{\et,Z}(X,J^\bullet(\mathcal{F})) := \textup{Cone}\left(\Gamma(X, J^\bullet(\mathcal{F})_\et) \to \Gamma(U, J^\bullet(\mathcal{F})_\et)  \right)[-1] \, .
\]

\begin{lem} \label{lem comparison etale}
\hfill \begin{enumerate} 
	\item Let $\mathcal{F}$ be a sheaf of abelian groups over $U_{\textup{\'et}}$. Then there is a canonical commutative diagram in the derived category of abelian groups, where the rows
are exact triangles:
\begin{changemargin}{-2cm}{5cm} 
\[\xymatrix{R\Gamma_{\et,c}(U,\mathcal{F}) \ar[r] \ar[d] & R\Gamma_{\et}(U, \mathcal{F}) \ar[r] \ar[d]^{\sim} & R\Gamma_{\et,Z}(X,j_! \mathcal{F})[1] \oplus \bigoplus_{v \in \Omega_{\R}} \widehat{R\Gamma_{\et}}(K_v,\mathcal{F}) \ar[r] \ar[d] &  R\Gamma_{\et,c}(U,\mathcal{F})[1] \ar[d] \\
R\Gamma_c(U, a^* \mathcal{F}) \ar[r] & R\Gamma(U, a^* \mathcal{F}) \ar[r] & R\Gamma(Z', i^* a^* \mathcal{F}) \oplus \bigoplus_{v \in \Omega_{\R}} \widehat{R\Gamma}(K_v, a^*\mathcal{F}) \ar[r] & R\Gamma_c(U, a^* \mathcal{F})[1] \, .
}
\]
\end{changemargin}

Besides, the complex $R\Gamma_{\et,Z}(X,j_! \mathcal{F})[1]$ 
is quasi-isomorphic to $\bigoplus_{v \in Z} R\Gamma_{\et}(K_v, {\mathcal F})$.

	\item Let $G$ be a smooth commutative group scheme over $U$.
Let $\underline{G}$ denote the fppf sheaf associated to $G$
and $G_\et := a_* \, \underline{G}$. Then there is a canonical commutative diagram in the derived category of abelian groups, where the rows
are exact triangles:
\begin{changemargin}{-2cm}{5cm} 
 \[\xymatrix{R\Gamma_{\et,c}(U,G_\et) \ar[r] \ar[d] & R\Gamma_{\et}(U, G_\et) \ar[r] \ar[d]^\sim & R\Gamma_{\et,Z}(X,j_! G_\et)[1] \oplus \bigoplus_{v \in \Omega_{\R}} \widehat{R\Gamma_{\et}}(K_v,G_\et) \ar[r] \ar[d] &  R\Gamma_{\et,c}(U,G_\et)[1] \ar[d]\\
R\Gamma_c(U, \underline{G}) \ar[r] & R\Gamma(U, \underline{G}) \ar[r] & R\Gamma(Z', i^* \underline{G}) \oplus \bigoplus_{v \in \Omega_{\R}} \widehat{R\Gamma}(K_v, \underline{G}) \ar[r] & R\Gamma_c(U, \underline{G})[1] \, .
}
\]
\end{changemargin}
Besides, the complex $R\Gamma_{\et,Z}(X,j_! G_\et)[1]$ 
is quasi-isomorphic to $\bigoplus_{v \in Z} R\Gamma_{\et}(K_v, G_\et)$.
\end{enumerate}
\end{lem}

\dem{}
\begin{enumerate}
	\item Set $J := J^\bullet(\mathcal{F})$. Since $j_! \mathcal{F} \to J_\et:=J^\bullet(\mathcal{F})_\et$ is an injective resolution, we get an injective resolution $\mathcal{F} = j^* j_! \mathcal{F} \to j^* J_\et$ in $U_{\et}$. The functor $\varepsilon^*$ is an exact functor that maps flasque \'etale sheaves to flasque fppf sheaves (see \cite[\href{http://stacks.math.columbia.edu/tag/0DDU}{Tag 0DDU}]{SP}), we get a flasque resolution $a^* \mathcal{F} \to I := \varepsilon^* j^* J_\et$. Let
${\widehat J}_v:={\widehat J}^{\bullet}({\mathcal F}_v)$; define
${\widehat I}_v={\widehat I}^{\bullet}((\varepsilon^* {\mathcal F})_v)$
(associated to the flasque resolution $I$ of 
$\varepsilon^* \mathcal{F}$) as in the beginning of section~\ref{sect1}. 

\smallskip

Consider now the following commutative diagram of complexes, where 
$\widetilde{\Gamma}_{\et,Z}(U,J$) and $\widetilde{\Gamma}_{\et,c}(U,J)$
are mapping cones defined such that the third and fourth rows are 
exact triangles:
\begin{changemargin}{-2cm}{5cm} 
\[\xymatrix{
\Gamma_{\et,c}(U,J_\et) \ar[r] & \Gamma_{\et}(U, j^*J_\et) \ar[r] & \Gamma_{\et,Z}(X,J_\et)[1] \oplus \bigoplus_{v \in \Omega_{\R}} \Gamma_{\et}(K_v,\widehat{J_\et}_v) \ar[r] &  \Gamma_{\et,c}(U, J_\et)[1] \\
\Gamma_{\et,c}(U,J) \ar[r] \ar[u]_{\varphi_c} \ar[d]^{d'} & \Gamma_{\et}(U, j^*J) \ar[r] \ar[u]_\varphi \ar[d]^= & \Gamma_{\et,Z}(X,J)[1] \oplus \bigoplus_{v \in \Omega_{\R}} \Gamma_{\et}(K_v,\widehat{J}_v) \ar[r] \ar[u]_{\varphi'} \ar[d]^d &  \Gamma_{\et,c}(U, J)[1] \ar[d] \ar[u] \\
\widetilde{\Gamma}_{\et,Z}(U,J) \ar[r] & \Gamma_{\et}(U, j^*J) \ar[r] & \bigoplus_{v \in Z} \Gamma_{\et,v}(\mathcal{O}_v,J)[1] \oplus \bigoplus_{v \in \Omega_{\R}} \Gamma_{\et}(K_v,\widehat{J}_v) \ar[r] &  \widetilde{\Gamma}_{\et,Z}(U, J)[1] \\
\widetilde{\Gamma}_{\et,c}(U,J) \ar[r] \ar[u]_{b'} \ar[d] & \Gamma_{\et}(U, j^*J) \ar[r] \ar[u]_= \ar[d]^{c} & \bigoplus_{v \in Z} \Gamma_{\et}(K_v,J) \oplus \bigoplus_{v \in \Omega_{\R}} \Gamma_{\et}(K_v,\widehat{J}_v) \ar[r] \ar[u]_b \ar[d] &  \widetilde{\Gamma}_{\et,c}(U, J)[1] \ar[u] \ar[d] \\
\Gamma_c(U, I) \ar[r] & \Gamma(U,I) \ar[r] & \Gamma(Z', i^* I) \oplus \bigoplus_{v \in \Omega_{\R}} \Gamma(K_v,\widehat{I}_v) \ar[r] & \Gamma_c(U, I)[1] \, .
}\]
\end{changemargin}
In this diagram, the rows are exact triangles (by definition for the last three rows, using the proof of Lemma 2.7 in \cite{GS} for the first ones). The maps $\varphi$, $\varphi'$ and $\varphi_c$ are quasi-isomorphisms by \cite[\href{http://stacks.math.columbia.edu/tag/0DDH}{Tag 0DDH}]{SP}. In addition, the maps $d$
and $b$ (hence also $d'$ and $b'$) are quasi-isomorphisms: for the map $d$, this is the excision property for \'etale cohomology (see \cite{MilEC}, Proposition III.1.27); for the map $b$, this is exactly \cite{MilADT}, Proposition II.1.1.(a). In addition, the  map $c$ is a quasi-isomorphism, using \cite[\href{http://stacks.math.columbia.edu/tag/0DDU}{Tag 0DDU}]{SP}. This proves the lemma.
\item Consider the following commutative diagram of exact triangles in the derived category
\begin{changemargin}{-2cm}{5cm} 
\[\xymatrix{
R\Gamma_{\et,c}(U, G_\et) \ar[r] \ar[d] & R\Gamma_{\et}(U, G_\et) \ar[r] \ar[d] & R\Gamma_{\et,Z}(X,j_! G_\et)[1] \oplus \bigoplus_{v \in \Omega_{\R}} \widehat{R\Gamma_{\et}}(K_v,G_\et) \ar[r] \ar[d] &  R\Gamma_{\et,c}(U,G_\et)[1] \ar[d]\\
R\Gamma_c(U, a^*  G_\et) \ar[r] \ar[d] & R\Gamma(U, a^*  G_\et) \ar[r] \ar[d] & R\Gamma(Z', i^* a^*  G_\et) \oplus \bigoplus_{v \in \Omega_{\R}} \widehat{R\Gamma}(K_v, a^* G_\et) \ar[r] \ar[d] & R\Gamma_c(U, a^*  G_\et)[1] \ar[d] \\
R\Gamma_c(U, \underline{G}) \ar[r] & R\Gamma(U, \underline{G}) \ar[r] & R\Gamma(Z', i^* \underline{G}) \oplus \bigoplus_{v \in \Omega_{\R}} \widehat{R\Gamma}(K_v, \underline{G}) \ar[r] & R\Gamma_c(U, \underline{G})[1] \, ,
}
\]
\end{changemargin}
where the vertical maps between the first two rows come from the first point of this Lemma, and the ones between the last two rows come from the adjunction morphism $a^* G_\et  = a^* a_* \underline{G} \to \underline{G}$ and from the functoriality of the triangle defining the complexes $R \Gamma_c(U,\cdot)$.
Now \cite{Brauer3}, Theorem 11.7, ensures that the composed vertical morphism $R\Gamma_\et(U,G_\et) \to R\Gamma(U,\underline{G})$ is an isomorphism. Whence the required result.
\end{enumerate}
\enddem

\dem{\ of Proposition \ref{prop supp compact}}
\begin{enumerate}
	\item This is immediate from the definitions, 
cf.~\cite{MilADT}, III, Proposition 0.4.a) and Remark 0.6.~b).

	\item The claim follows from the definitions,
from the exactness of the functors $i^*$, $a^v_*$
and $D^\bullet(\cdot)$ at the beginning of section \ref{sect1}, and from the exactness of the cone functor on the category of complexes of abelian groups (see also \cite{MilADT}, III, Proposition 0.4.b) and Remark 0.6.~b)).

	\item As in the proof of  \cite{MilADT}, III, Proposition 0.4.c), let $I^\bullet(\mathcal{F})$ be an injective resolution of $\mathcal{F}$. 
In the number field case, {\it the piece of notation $\Gamma(\widehat{K_v}, 
I^\bullet(\mathcal{F}))$ will stand for 
$\Gamma(K_v, {\widehat I}^\bullet(\mathcal{F}_v))$ when $v$ is a real place 
of $K$}, where ${\widehat I}^\bullet(\mathcal{F}_v)$ is the modified 
resolution constructed in the beginning of section~\ref{sect1}.

\smallskip

Consider the following commutative diagram of complexes in the category
of bounded below complexes of abelian groups:
\[
\xymatrix{
\Gamma(U, I^\bullet(\mathcal{F})) \ar[r]^-{\alpha} \ar[d]^-f & \bigoplus_{v \notin U} \Gamma(\widehat{K_v}, I^\bullet(\mathcal{F})) \oplus \bigoplus_{v \in U \setminus V} \Gamma(\widehat{\mathcal{O}_v}, I^\bullet(\mathcal{F})) \ar[d]^-{(\id, g)} \ar[r]^-{\pi_{\mathcal{O}}} & \bigoplus_{v \notin U} \Gamma(\widehat{K_v}, I^\bullet(\mathcal{F})) \\
\Gamma(V,I^\bullet(\mathcal{F})) \ar[r]^-{\alpha'} & \bigoplus_{v \notin U} \Gamma(\widehat{K_v}, I^\bullet(\mathcal{F})) \oplus \bigoplus_{v \in U \setminus V} \Gamma(\widehat{K_v}, I^\bullet(\mathcal{F})) \ar[ru]_-{\pi_K} \, ,
}
\]
where the maps are the natural ones.

Functoriality of the mapping cone in the category of complexes gives morphisms
$$\Gamma_{U \setminus V}(U, I^\bullet(\mathcal{F})) \to \bigoplus_{v \in U \setminus V} \Gamma_v(\mathcal{O}_v, I^\bullet(\mathcal{F})) \to \bigoplus_{v \in U \setminus V} \Gamma_v(\widehat{\mathcal{O}_v}, I^\bullet(\mathcal{F})) \, ,$$
where $$\Gamma_{U \setminus V}(U, I^\bullet(\mathcal{F})) := C(f)[-1], \quad 
\Gamma_v(\mathcal{O}_v, I^\bullet(\mathcal{F})) := \Gamma_{\mathfrak{m}_v}(\mathcal{O}_v, I^\bullet(\mathcal{F}))$$ and $$\Gamma_v(\widehat{\mathcal{O}_v}, I^\bullet(\mathcal{F})) := \Gamma_{\mathfrak{m}_v}(\widehat{\mathcal{O}_v}, I^\bullet(\mathcal{F})).$$

The excision property (Lemma \ref{lem excision}) implies that the first morphism $\Gamma_{U \setminus V}(U, I^\bullet(\mathcal{F})) \to \bigoplus_{v \in U \setminus V} \Gamma_v(\mathcal{O}_v, I^\bullet(\mathcal{F}))$ is a quasi-isomorphism.

Since for all $v \in X$, the ring $\mathcal{O}_v$ is an excellent henselian discrete valuation ring, Lemma \ref{lem artin} ensures that the second map 
$$\bigoplus_{v \in U \setminus V} \Gamma_v(\mathcal{O}_v, I^\bullet(\mathcal{F})) \to \bigoplus_{v \in U \setminus V} \Gamma_v(\widehat{\mathcal{O}_v}, I^\bullet(\mathcal{F}))$$
is a quasi-isomorphism. Therefore, the natural morphism $C(f) \to C(g)$ is a quasi-isomorphism.

Apply now Lemma \ref{lem complex} to get a 
commutative diagram in the derived category of abelian groups:
\begin{changemargin}{-4cm}{5cm} 
\begin{equation} \label{bigdiag}
\begin{gathered}
\resizebox{20cm}{!}{
\xymatrix{
\left(\bigoplus_{v \notin V} \Gamma(\widehat{K_v}, I^\bullet(\mathcal{F}))\right)[-1] \ar[d] & \left(\bigoplus_{v \notin U} \Gamma(\widehat{K_v}, I^\bullet(\mathcal{F}))\right)[-1] \ar[l]_-{i_K} \ar[d] & \bigoplus_{v \notin U} \Gamma(\widehat{K_v}, I^\bullet(\mathcal{F})) \oplus \bigoplus_{v \in U \setminus V} \Gamma(\widehat{\mathcal{O}_v}, I^\bullet(\mathcal{F})) \ar[r]^-{(\id, g)} & \bigoplus_{v \notin V} \Gamma(\widehat{K_v}, I^\bullet(\mathcal{F})) \ar[d] \\
\Gamma_c(V,I^\bullet(\mathcal{F})) \ar[r] \ar[d] & \Gamma_c(U,I^\bullet(\mathcal{F})) \ar[r] \ar[d] & \bigoplus_{v \in U \setminus V} \Gamma(\widehat{\mathcal{O}_v}, I^\bullet(\mathcal{F}))  \ar[u]^{i'_{\mathcal{O}}} \ar[r] & \Gamma_c(V,I^\bullet(\mathcal{F}))[1] \\
\Gamma(V,I^\bullet(\mathcal{F})) \ar[d]^-{\alpha'} & \Gamma(U,I^\bullet(\mathcal{F})) \ar[d]^-{\pi_{\mathcal{O}} \circ \alpha} \ar[l]_{f} \ar[r]^-{\alpha} & \bigoplus_{v \notin U} \Gamma(\widehat{K_v}, I^\bullet(\mathcal{F})) \oplus \bigoplus_{v \in U \setminus V} \Gamma(\widehat{\mathcal{O}_v}, I^\bullet(\mathcal{F}))  \ar[u]^{\pi'_{\mathcal{O}}} & \\
\bigoplus_{v \notin V} \Gamma(\widehat{K_v}, I^\bullet(\mathcal{F})) \ar[r]^-{\pi_K} & \bigoplus_{v \notin U} \Gamma(\widehat{K_v}, I^\bullet(\mathcal{F})) \, , & & 
}
}
\end{gathered}
\end{equation}
\end{changemargin}
where the second row and the first two columns are exact triangles.

Now the cohomology of this diagram gives the following canonical 
commutative diagram, with an exact second row 
(and the two first columns exact):
\begin{changemargin}{-4cm}{4cm} 
 \[
\resizebox{20cm}{!}{
\xymatrix{
& \bigoplus_{v \notin V} H^{r-1}(\widehat{K_v},\mathcal{F}) \ar[d] & \bigoplus_{v \notin U} H^{r-1}(\widehat{K_v}, \mathcal{F}) \ar[l] \ar[d] & \bigoplus_{v \notin U} H^r(\widehat{K_v}, \mathcal{F}) \oplus \bigoplus_{v \in U \setminus V}  H^r(\widehat{\mathcal{O}_v}, \mathcal{F})  \ar[r] & \bigoplus_{v \notin V} H^r(\widehat{K_v}, \mathcal{F}) \ar[d] & \\
\dots \ar[r]& H^r_c(V, \mathcal{F}) \ar[r] \ar[d] & H^r_c(U, \mathcal{F}) \ar[r] \ar[d] & \bigoplus_{v \in U \setminus V} H^r(\widehat{\mathcal{O}_v}, \mathcal{F}) \ar[u] \ar[r] & H^{r+1}_c(V,\mathcal{F}) \ar[r] & \dots\\
& H^r(V,\mathcal{F}) \ar[d]& H^r(U,\mathcal{F}) \ar[d] \ar[l]_-{\textup{Res}} \ar[r] & \bigoplus_{v \notin U} H^r(\widehat{K_v}, \mathcal{F}) \oplus \bigoplus_{v \in U \setminus V}   H^r(\widehat{\mathcal{O}_v}, \mathcal{F})  \ar[u] & & \\
& \bigoplus_{v \notin V} H^r(\widehat{K_v},\mathcal{F}) \ar[r] & \bigoplus_{v \notin U} H^r(\widehat{K_v}, \mathcal{F}) \, , & & &
}
}
\]
\end{changemargin}

\smallskip

which proves the required exactness and commutativity.

	\item Lemma \ref{lem comparison etale} gives a commutative diagram
 with exact rows:
\[
\xymatrix{
H^{r-1}_{\et}(U, \mathcal{F}) \ar[r] \ar[d]^{\sim} & \bigoplus_{v \notin U} H^{r-1}_\et(K_v, \mathcal{F}) \ar[r] \ar[d] & H^r_{\et,c}(U, \mathcal{F}) \ar[r] \ar[d] & H^r_{\et}(U, \mathcal{F}) \ar[r] \ar[d]^{\sim} & \bigoplus_{v \notin U} H^{r}_\et(K_v, \mathcal{F}) \ar[d] \\
H^{r-1}(U, \mathcal{F}) \ar[r] & \bigoplus_{v \notin U} 
H^{r-1}(\widehat{K_v}, \mathcal{F}) \ar[r] & H^r_{c}(U, \mathcal{F}) 
\ar[r] & H^r(U, \mathcal{F}) \ar[r] & 
\bigoplus_{v \notin U} H^{r}(\widehat{K_v}, \mathcal{F})  .
}
\]
Here $H_{\et}$ stands for \'etale cohomology (modified over 
$K_v$ for $v$ real) and $H_{\et,c}$ for (modified) \'etale
cohomology with compact support (as defined in \cite{GS}, \S 2, or before Lemma \ref{lem comparison etale}; recall also
that in the number field case, the piece of notation $v \not \in U$ means 
that we consider the places corresponding to closed points of $\spec (\calo_K) 
\setminus U$ {\it and} the real places).

By \cite{GGMB}, Lemma 3.5.3, and \cite{MilEC} III.3, we have 
$$H^r_\et(K_v, \mathcal{F})\cong H^r_\et(\widehat K_v, \mathcal{F}) \xrightarrow{\sim} H^r(\widehat{K_v}, 
\mathcal{F})$$ for all $r \geq 1$ (resp.~for all integers $r$ if 
$\mathcal{F}_K$ is finite; indeed $K_v$ and $\widehat{K_v}$ have the same 
absolute Galois group via \cite{bkiac}, \S 8, Corollary 4 to Theorem 2 and 
\cite{riben}, section F, Cor.~2 to Th.~2)
and all places $v$ of $K$.  
Therefore the five-lemma gives the result. \enddem
\end{enumerate}

\begin{rema}
{\rm The definition of fppf compact support cohomology and its related
properties are specific to schemes of dimension $1$. To the best of our
knowledge, there is no good analogue in higher dimension, unlike what
happens for \'etale cohomology.
}
\end{rema}

We will need the following complement to Proposition~\ref{prop supp compact}:
\begin{prop} \label{bonus}
Let ${\cal F}$ be a flat affine commutative group scheme 
of finite type over $U$. Let $V \subset U$ be a non empty 
open subset. Then there is a long exact sequence
\begin{equation} \label {firstlong}
\dots \to \bigoplus_{v \in U \setminus V} H^r_v(\widehat{\mathcal{O}_v}, \mathcal{F}) \to
H^r(U, \mathcal{F}) \to
H^r(V,\mathcal{F}) \to \bigoplus_{v \in U \setminus V} H^{r+1}_v(\widehat{\mathcal{O}_v}, \mathcal{F})
\to \dots
\end{equation}
\end{prop}

\dem{} 
The map $\bigoplus_{v \in U \setminus V}
H^r_v(\widehat{\mathcal{O}_v}, \mathcal{F}) \to H^r(U, \mathcal{F})$ is
given by the
identification of the first group with $H^r_Z(U,\mathcal{F})$, where $Z=U \setminus V$, via Lemma~\ref{lem artin} and Lemma~\ref{lem excision}.
By the localization exact sequence (\cite{MilADT}, Prop.~III.0.3.~c),
this identification yields the required long exact sequence. 
\enddem

\section{Topology on cohomology groups with compact support} \label{section topo}

With the previous notation, let us define a natural topology on the 
groups $H^*_c(U,N)$, where $N$ is a finite flat commutative $U$-group scheme.
Th.~\ref{thm AV} actually immediately implies that 
$H^2_c(U,N)$ is profinite, but this duality theorem will not be used 
in this paragraph. The ``a priori'' approach we adopt in this section 
answers a question raised by Milne (\cite{MilADT}, Problem III.8.8.).

We restrict ourselves to the function field case, because when $K$ is a number
field the groups involved are finite (cf.~\cite{MilADT}, Th.~III.3.2; see also 
section~\ref{sect3} of this article).
Recall that as usual (cf.~for example \cite{MilADT}, \S III.8),
the groups $H^*(U, N)$ are endowed with the discrete topology. Our first goal in this section is to define a natural topology on the 
groups $H^*_c(U,N)$. 

\smallskip

Given an exact sequence of abelian groups
$$0 \to A \to B \to C \to 0 \, ,$$
such that $A$ is a topological group, there exists a unique 
topology on $B$ such that $B$ is a topological group, $A$ is an open subgroup of $B$, and $C$ is discrete when endowed with the quotient topology. Indeed, the topology on $B$ is generated by the subsets $b + U$, where $b \in B$ and $U$ is an open subset of $A$. In addition, given another abelian group $B'$ with a subgroup $A' \subset B'$ that is a topological group, and a commutative diagram of abelian groups
\[
\xymatrix{
A \ar@{^{(}->}[r] \ar[d]^-f & B \ar[d]^-g \\
A' \ar@{^{(}->}[r] & B' \, ,
}
\]
then $f$ is continuous if and only if $g$ is continuous, for the aforementioned topologies. And $f$ is open if and only if $g$ is.

We can therefore topologize the groups
$H^i_c(U, N)$ for $i \neq 2$, using the exact sequence
(see Proposition \ref{prop supp compact}, 1.)
$$ \bigoplus_{v \not \in U} H^{i-1}(\widehat{K_v},N) \to H^i_c(U, N) \to
H^i(U,N).$$ Since the groups $H^{i-1}(\widehat{K_v},N)$ are finite
for $i \neq 2$ (\cite{MilADT}, \S III.6)
and $H^i(U,N)$ is discrete, all groups $H^i_c(U, N)$ are discrete if $i \neq 2$.

\smallskip 

Let us now focus on the case $i=2$.
Consider the exact sequence (Proposition \ref{prop supp compact}, 1.)
\begin{equation} \label{bonusuite}
H^1(U, N) \to \bigoplus_{v \notin U} H^1(\widehat{K}_v,N) \to H^2_c(U, N) \to H^2(U, N) \, .
\end{equation}
and for $i=1,2$, set 
$$D^i(U,N)={\rm Im}\, [H^i_c(U,N) \to H^i(U,N)]={\rm Ker} \,
[H^i(U,N) \to \bigoplus_{v \not \in U} H^i(\widehat{K_v},N)].$$
By Proposition~\ref{prop supp compact}, 1., there is an exact sequence
\begin{equation} \label {didef}
\bigoplus_{v \not \in U} H^{i-1}(\widehat{K_v},N) \to H^i_c(U,N) \to
D^i(U,N) \to 0.
\end{equation}

The following result has been proved by \v{C}esnavi\v{c}ius (\cite{cesnalms}, 
Th.~2.9).\footnote {Proposition~2.3 of loc.~cit.~uses the fppf duality 
Theorem~\ref{thm AV}, but this proposition is actually not needed to prove 
Theorem~\ref{strict} because a discrete subgroup of a Hausdorff topological
group is automatically closed by \cite{tg}, \S 2, Prop.~5.}

\begin{theo}[\v{C}esnavi\v{c}ius] \label{strict}
The map $H^1(U,N) \to \bigoplus_{v \notin U} 
H^1(\widehat{K}_v,N)$ is a strict morphism of topological groups,
that is: the image of $H^1(U,N)$ is a discrete subgroup of 
$\bigoplus_{v \notin U} H^1(\widehat{K}_v,N)$. Besides, the 
group $D^1(U,N)$ is finite.
\end{theo}

\begin{cor} \label{finih1c}
The group $H^1_c(U,N)$ is finite. 
\end{cor}

\dem{} The group $\bigoplus_{v \notin U} H^0(\widehat{K_v},N)$ is finite 
($N$ being a finite $U$-group scheme). Thus the finiteness 
of $H^1_c(U,N)$ is equivalent to the finiteness of $D^1(U,N)$ 
by (\ref{didef}).
\enddem

Put the quotient topology on 
$(\bigoplus_{v \not \in U} H^1(\widehat{K}_v,N))/{\rm Im} \, 
H^1(U,N)$. Using Th.~\ref{strict},
the previous facts define a natural topology on $H^2_c(U, N)$, so that morphisms in the exact sequence (\ref{bonusuite}) are continuous (and even 
strict). 
This topology makes $H^2_c(U, N)$ a Hausdorff and locally compact group
(cf.~\cite{tg}, \S 2, Prop.~18, a).

\smallskip

To say more about the topology of $H^2_c(U, N)$, we need a lemma:

\begin{lem} \label{continuuv}
\hfill \begin{enumerate} \item Let $r: N \to N'$ be a morphism of finite flat commutative 
$U$-group schemes. Then 
the corresponding map $s : H^2_c(U,N) \to H^2_c(U,N')$ is continuous. If 
we assume further that $r$ is surjective, then $s$ is open. If 
$$0 \to N' \to N \to N'' \to 0$$
is an exact sequence of finite flat commutative
$U$-group schemes, then the connecting map $H^2_c(U,N'') \to H^3_c(U,N')$
is continuous.

	\item Let $V \subset U$ be a non empty open subset. Then 
the natural map $u: H^2_c(V,N) \to H^2_c(U, N)$ is continuous.
\end{enumerate}
\end{lem}

\dem{} \hfill \begin{enumerate}
	\item By definition of the topology on the groups $H^2_c$, it is 
sufficient to prove that for $v \not \in U$, 
the map $H^1(\widehat{K_v},N) \to H^1(\widehat{K_v},N')$ is continuous 
(resp.~open if $r$ is surjective). 
Continuity follows from \cite{Ces}, Prop.~4.2 and the openness statement
from loc.~cit., Prop~4.3 d). Similarly, the last assertion follows
from the continuity of the connecting map $H^1(K_v,N'') \to H^2(K_v,N')$
(loc.~cit.,  Prop.~4.2).

	\item Since (by definition of the topology) the image $I$ of 
$A:=\bigoplus_{v \not \in V} H^1(\widehat{K_v},N)$ is an open subgroup 
of $H^2_c(V,N)$, it is sufficient to show that the restriction 
of $u$ to $I$ is continuous. As $I$ is equipped with the quotient 
topology (induced by the topology of $A$), this is equivalent to showing that 
the natural map $s: A \to H^2_c(U, N)$
is continuous. Now we observe that  $A$ is the direct sum of 
$A_1:= \bigoplus_{v \not \in U} H^1(\widehat{K_v},N)$ and 
$A_2:=\bigoplus_{v \in U \setminus V} H^1(\widehat{K_v},N)$.
The restriction of $s$ to $A_1$ is continuous by the commutative diagram
of Prop.~\ref{prop supp compact}, 3. Therefore it only remains to show that 
the restriction $s_2$ of $s$ to $A_2$ is continuous. By loc.~cit., 
the restriction of $s_2$ to $\bigoplus_{v \in U \setminus V} 
H^1(\widehat{\mathcal{O}_v},N)$ is zero. Since $\bigoplus_{v \in U \setminus V} 
H^1(\widehat{\mathcal{O}_v},N)$ is an open subgroup of $\bigoplus_{v \in U \setminus V} 
H^1(\widehat{K_v},N)$ (\cite{Ces}, Prop.~3.10), the result follows.
\end{enumerate}
\enddem 

Recall also the following (probably well-known) lemma: 

\begin{lem} \label{strictprof}
Let $f : A \to B$ be a continuous morphism of topological groups, with 
$B$ Hausdorff. 
\begin{enumerate}
	\item Assume that $A$ is profinite. Then $f$ is strict.
	\item Assume that $f$ is injective and $A$ is compact. Then $f$ is strict. 
	\item Let 
$$0 \to A \stackrel{i}{\to} B \stackrel{\pi}{\to} C \to 0$$
be an exact sequence of topological groups with $i$ strict and $\pi$
continuous. If $A$ and $C$ are completely disconnected, then so is $B$.
\end{enumerate}
\end{lem}

\dem{} \hfill \begin{enumerate}
	\item Since $f$ is continuous and $B$ Hausdorff, the image of $f$ is a 
compact subspace of $B$, so we can assume that $B$ is compact and $f$ 
is onto. The topology of $A$ has a basis consisting of open subgroups,
so it is sufficient to show that the image of such a subgroup $U$ 
is open. As $U$ is closed (hence compact) and of finite index in $A$,
its image $f(U)$ is also compact and of finite index in $B$, hence it is 
an open subgroup of $B$.

	\item Since $A$ is compact and $B$ is Hausdorff, we get that $i$ is a 
closed map (because the image of a compact subspace of $A$ is compact),
hence it induces a homeomorphism from $i$ onto the subspace 
$i(A) \subset B$. This means that $i$ is strict.

	\item Let $D$ be a connected subset of $B$. Then $\pi(D)$ is connected, 
hence is a singleton. Thus, by translating, one can assume that $D \subset i(A)$; as $i$ is strict, the subset 
$i^{-1}(D) \subset A$ is connected,
so it is reduced to a point, hence $D$ is a singleton. This proves the statement.
\end{enumerate}
\enddem

\begin{prop} \label{profiniteh2}
For every integer $i$ with $0 \leq i \leq 3$, 
the topology on $H^i_c(U, N)$ is profinite. 
\end{prop}

\dem{} 
The only non trivial case is $i=2$.
We first observe that if there is an exact sequence of finite
flat commutative $U$-group schemes 
$$0 \to N' \to N \to N'' \to 0,$$
then it is sufficient to prove that $H^2_c(U,N')$ and 
$H^2_c(U,N'')$ are profinite to get the same result for $H^2_c(U,N)$. 
Indeed by Proposition~\ref{prop supp compact}, 3., there is an exact sequence
$$ H^1_c(U,N'') \to H^2_c(U,N') \to H^2_c(U,N) \to H^2_c(U,N'').$$
The group $H^1_c(U,N'')$ is finite by Corollary~\ref{finih1c}; besides, 
all maps are continuous and the map $H^2_c(U,N) \to H^2_c(U,N'')$ is open 
(in particular it is strict, and 
its image is profinite as soon as $H^2_c(U,N'')$ is) 
by  Lemma~\ref{continuuv}, 1. 
Therefore if $H^2_c(U,N')$ and $H^2_c(U,N'')$ are profinite, 
then $H^2_c(U,N)$ is profinite as an extension 
$$0 \stackrel{i}{\to} A \to H^2_c(U,N) \stackrel{\pi}{\to} B \to 0$$
of two profinite groups $A$, $B$ such that the map $\pi$ is open
(the map $i$ is strict by Lemma~\ref{strictprof}, 2.;
the group $H^2_c(U,N)$ is completely disconnected by
Lemma~\ref{strictprof} 3., and 
its compactness follows from the fact that $\pi$ is a proper map 
by \cite{tg}, \S 4, Cor.~2 to Prop.~2).

\smallskip 

This being said, note now that Proposition \ref{prop supp compact}, 4.~implies the result when the order of 
$N$ is prime to $p$ by \cite{MilADT}, Corollary II.3.3 (in this case 
$H^2_c(U,N)$ is even finite). One can therefore assume by devissage 
that the order of $N$ is a power of $p$.
The generic fiber $N_K$ of $N$ is a finite commutative group scheme over $K$. 
By \cite{demgab}, IV, \S 3.5, $N_K$
 admits a composition series whose quotients are \'etale (with a dual of height one), local (of height one) with \'etale dual, or $\alpha_p$. The schematic closure in $N$ of this composition series provides a composition series defined over $U$.
Thus, using the same devissage argument as above, one reduces to the case where 
the generic fiber $N_K$ or its dual $N_K^D$ has height one. 

Proposition III.B.4 and Corollary III.B.5
in \cite{MilADT} now imply that there exists a non empty open subset $V \subset U$ such that $N_{\vert V}$ extends to a finite flat commutative 
group scheme $\widetilde{N}$ over the proper $k$-curve $X$. 

Then Proposition \ref{prop supp compact}, 3.~gives an exact sequence
\begin{equation} \label{ex seq 3}
H^1_c(X, \widetilde{N}) \to \bigoplus_{v \in X \setminus V} H^1(\widehat{\mathcal{O}}_v, \widetilde{N}) \to H^2_c(V,N) \to H^2_c(X, \widetilde{N}) \, 
\end{equation}
and since we are in the function field case with $X$ proper over $k$, 
we have $H^i_c(X, \widetilde{N})=H^i(X,\widetilde{N})$ for every 
positive integer $i$.

By Proposition \ref{prop supp compact}, 3., the map 
$\bigoplus_{v \in X \setminus V} H^1(\widehat{\mathcal{O}}_v, \widetilde{N}) 
\to H^2_c(V,N)$ factors through 
$\bigoplus_{v \in X \setminus V} H^1(\widehat{K}_v,N)$, hence it is continuous. 
By Lemma~\ref{continuuv}, 
all maps in \eqref{ex seq 3} are continuous. In addition, the groups 
$H^1_c(X, \widetilde{N})=H^1(X, \widetilde{N})$ and 
$H^2_c(X, \widetilde{N})=H^2(X, \widetilde{N})$ are finite by 
\cite{MilADT}, Lemma~III.8.9. Besides,
$\bigoplus_{v \in X \setminus V} H^1(\widehat{\mathcal{O}}_v, \widetilde{N})$ is profinite by loc.~cit., \S III.7; hence $H^2_c(V,N)$ is profinite as 
an extension
(the maps being strict by Lemma~\ref{strictprof}, 2.)
of a finite group by a profinite group.

\smallskip

Since $H^2(\widehat{\mathcal{O}}_v, N)=0$ for every $v \in U$
(\cite{MilADT}, \S III.7), Prop.~\ref{prop supp compact}, 3.~
gives an exact sequence of groups
$$ \bigoplus_{v \in U \setminus V} H^1(\widehat{\mathcal{O}}_v, N) \to H^2_c(V,N) \to H^2_c(U, N) \to 0, $$
which implies that $H^2_c(U, N)$ is profinite, the map 
$H^2_c(V,N) \to H^2_c(U, N)$ being continuous by Lemma~\ref{continuuv} 2.,
hence strict by Lemma~\ref{strictprof} 1.,
because $H^2_c(V,N)$ is profinite and $H^2_c(U, N)$ is Hausdorff. 
\enddem

The following statement will be useful in the next section:

\begin{prop} \label{rem topo}
Assume that $\mathcal F=N$, $\mathcal F'=N'$ and $\mathcal F''=N''$ are 
finite and flat commutative 
group schemes over $U$.
Then all the maps in Proposition \ref{prop supp compact} are strict 
(in particular continuous). 
\end{prop}

\dem{} For the maps in assertion 1.~of Prop.~\ref{prop supp compact}, 
this follows from the definition of the topology and Th.~\ref{strict}. 

\smallskip

Let us consider the maps in assertion 2.
The finiteness of the $H^1_c$ groups (Cor.~\ref{finih1c}) implies 
that it only remains to deal with the maps between $H^2_c$'s and the 
connecting map $H^2_c(U,{\mathcal F}'') \to H^3_c(U, {\mathcal F}')$.
All these maps are continuous by Lemma~\ref{continuuv}, hence
strict by Lemma~\ref{strictprof} 1. and Prop.~\ref{profiniteh2}.

\smallskip

Finally, it has already been proven (cf.~proof of Prop.~\ref{profiniteh2})
that the maps in the exact sequence of assertion 3.~are continuous. 
They are strict via Lemma~\ref{strictprof} 1.
because $H^1_c(U,{\mathcal F})$ is finite, 
$H^2_c(U,{\mathcal F})$ (resp.~$\bigoplus_{v \in U \setminus V} H^1(\widehat{\calo_v}, 
{\mathcal F})$) is profinite, and the other groups are discrete.
\enddem

\section{Proof of Theorem \ref{thm AV} in the function field 
case} \label{sect2}

In this section $K$ is the function field of a projective, smooth 
and geometrically integral curve $X$ defined over a finite field $k$ of 
characteristic $p$.
The proof follows the same lines as the proof of \cite{MilADT}, 
Theorem III.8.2, 
replacing Proposition III.0.4 in \cite{MilADT} by 
Proposition \ref{prop supp compact} and 
using the results of section~\ref{sect1}.

\smallskip

For every non empty open subset $V \subset U$, the natural map
$H^3_c(V,\G) \stackrel{s}{\to} H^3_c(U,\G)$ is an isomorphism, and the trace map
identifies $H^3_c(U,\G)$ with $\Q/\Z$
(this identification being compatible with $s$).
Indeed since $\G$ is a smooth group 
scheme we can apply Prop~\ref{prop supp compact}, 4.~
and \cite{MilADT}, \S II.3.

For a fppf sheaf ${\mathcal F}$ 
on $U$, let us first define the pairing of abelian groups
$$H^{3-r}_c(U, \mathcal{F}) \times H^{r}(U, \mathcal{F}^D) \to H^3_c(U, \G) \cong \Q / \Z \, .$$

Since the cohomology groups with compact support are defined via a mapping cone construction, we need to construct this pairing carefully at the level of complexes in order to be able to prove the compatibilities that follow (see Lemmas \ref{lem cup cob} and \ref{lem commutative} for instance).

\begin{lem} \label{lem cup}
Let $A$ and $B$ be two fppf sheaves of abelian groups 
on $U$. Then there exists a canonical pairing in the derived category of abelian groups:
$$R\Gamma_c(U,A) \otimes^{\L} R\Gamma(U,B) \to R\Gamma_c(U,A \otimes B) \, .$$
Moreover, this pairing is functorial in $A$ and $B$.
\end{lem}

\dem{}
For any complex $C$ of fppf sheaves, let $G(C)$ denote the Godement resolution of $C$ (see for instance \cite{SGA4}, XVII, 4.2.9; Godement resolutions exist on the big fppf site because this site has enough points, see Remark~1.6.~of \cite{gabkelly} or \cite[\href{http://stacks.math.columbia.edu/tag/06VX}{Tag 06VX}]{SP}). 

Then there is a commutative diagram of complexes of sheaves
(see \cite{God}, II.6.6 or \cite{friedsus}, Appendix A)
\[
\xymatrix{
A \otimes B \ar[rd] \ar[d] & \\
\textup{Tot}(G(A) \otimes G(B)) \ar[r] & G(A \otimes B) \, .
}
\]
The horizontal morphism induces a morphism of complexes of abelian groups
$$\textup{Tot}(\Gamma(U,G(A)) \otimes \Gamma(U,G(B))) \to \Gamma(U,G(A \otimes B)) \, $$
hence a canonical morphism in the derived category of abelian groups
$$\Gamma(U,G(A)) \otimes^{\L} \Gamma(U,G(B)) \to \Gamma(U,G(A \otimes B)) \, .$$

Considering the local versions of the previous pairings, one gets a commutative diagram of complexes of abelian groups 
\[
\xymatrix{
\textup{Tot}(\Gamma(U,G(A))\otimes \Gamma(U,G(B))) \ar[r] \ar[d] & \Gamma(U,G(A \otimes B)) \ar[d] \\
\prod_{v \notin U} \textup{Tot}(\Gamma(\widehat{K}_v, G(A)) \otimes \Gamma(U,G(B))) \ar[r] & \prod_{v \notin U} \Gamma(\widehat{K_v}, G(A \otimes B)) \, ,
}
\]
and functoriality of cones gives a canonical morphism of complexes
(via Proposition~\ref{tensorcomplex} in the Appendix)
\begin{equation} \label{cup complexes} 
\textup{Tot}(\Gamma_c(U, G(A)) \otimes \Gamma(U,G(B))) \to \Gamma_c(U,G(A \otimes B)) \, .
\end{equation}

Since Godement resolutions are acyclic (see \cite{SGA4}, XVII, Proposition 4.2.3), we know that $R\Gamma(U,C) \cong \Gamma(U,G(C))$ in the derived category, for any fppf sheaf $C$.  
Hence the pairing \eqref{cup complexes} gives the required morphism in the derived category
$$R\Gamma_c(U,A) \otimes^{\L} R\Gamma(U,B) \to R\Gamma_c(U,A \otimes B) \, .$$

The functoriality of Godement resolutions implies the functoriality of the pairing in $A$ and $B$. 
\enddem

Using Lemma \ref{lem cup}, \cite[\href{http://stacks.math.columbia.edu/tag/068G}{Tag 068G}]{SP}
gives a natural pairing
$$H^r_c(U,A) \times H^s(U,B) \to H^{r+s}_c(U, A \otimes B) \, ,$$
whence we deduce the required canonical pairings, for any sheaf $\mathcal{F}$ on $(\textup{Sch}/U)_{\textup{fppf}}$
\begin{equation} \label{cup}
H^r_c(U,{\mathcal F}) \times H^s(U,{\mathcal F}^D) \to H^{r+s}_c(U, \G) \, ,
\end{equation}
using the canonical map ${\mathcal F} \otimes \mathcal{F}^D = \mathcal{F} \otimes \underline{\Hom}(\mathcal{F}, \G) \to \G$.

\smallskip

Let us describe explicitely the pairing above: the map
\[
\textup{Tot}(\Gamma_c(U, G(A)) \otimes \Gamma(U,G(B))) \to \Gamma_c(U,G(A \otimes B))
\]
is given by maps 
\begin{changemargin}{-2cm}{5cm} 
\begin{align*}
\left(\prod_{v \notin U} \Gamma(\widehat{K}_v,G_{r-1}(A)) \oplus \Gamma(U, G_r(A))\right) \otimes \Gamma(U,G_s(B)) & \to \prod_{v \notin U} \Gamma(\widehat{K}_v,G_{r+s-1}(A \otimes B)) \oplus \Gamma(U, G_{r+s}(A \otimes B)) \\
(a_{r-1},a_r) \otimes b_s & \mapsto (a_{r-1} \cup \beta(b_s), a_{r} \cup b_s) \, ,
\end{align*}
\end{changemargin}
where the maps denoted by $\cup$ are the natural pairings, and $\beta : \Gamma(U,G_s(B)) \to \prod_{v \notin U} \Gamma(\widehat{K}_v,G_{s}(B))$ is the localization map.

In the following, we will need an alternative version of the above pairing: with the same notation as above, one defines a pairing in the derived category
\[
R\Gamma(U,A) \otimes^{\L} R\Gamma_c(U,B) \to R\Gamma_c(U,A \otimes B) \, .
\]
The definition is similar to the one in Lemma \ref{lem cup}: the commutative diagram of complexes 
\[
\xymatrix{
\textup{Tot}(\Gamma(U,G(A))\otimes \Gamma(U,G(B))) \ar[r] \ar[d] & \Gamma(U,G(A \otimes B)) \ar[d] \\
\prod_{v \notin U} \textup{Tot}(\Gamma(U, G(A)) \otimes \Gamma(\widehat{K}_v,G(B))) \ar[r] & \prod_{v \notin U} \Gamma(\widehat{K_v}, G(A \otimes B)) \, ,
}
\]
and Proposition \ref{tensorcomplex} in the Appendix gives a morphism of complexes
\begin{equation} \label{cup complexes right} 
\textup{Tot}(\Gamma(U, G(A)) \otimes \Gamma_c(U,G(B))) \to \Gamma_c(U,G(A \otimes B)) \, .
\end{equation}

Taking into account the signs in Proposition \ref{tensorcomplex}, one can describe the pairing $\textup{Tot}(\Gamma(U,G(A)) \otimes \Gamma_c(U,G(B))) \to \Gamma_c(U,G(A \otimes B))$ explicitely as follows:
\begin{changemargin}{-2cm}{5cm} 
\begin{align*}
\Gamma(U, G_r(A)) \otimes \left(\prod_{v \notin U} \Gamma(\widehat{K}_v,G_{s-1}(B)) \oplus \Gamma(U,G_s(B))\right) & \to \prod_{v \notin U} \Gamma(\widehat{K}_v,G_{r+s-1}(A \otimes B)) \oplus \Gamma(U, G_{r+s}(A \otimes B)) \\
 a_r \otimes (b_{s-1},b_{s}) & \mapsto ((-1)^r \alpha(a_r) \cup b_{s-1}, a_{r} \cup b_{s}) \, ,
\end{align*}
\end{changemargin}
where $\alpha : \Gamma(U,G_r(A)) \to \prod_{v \notin U} \Gamma(\widehat{K}_v,G_{r}(A))$ is the localization map.

We now compare the two pairings defined above:
\begin{lem} \label{lem exchange left right}
The following diagram of complexes
\[
\xymatrix{
\textup{Tot}(\Gamma_c(U,G(A)) \otimes \Gamma(U,G(B))) \ar[r] & \Gamma_c(U,G(A \otimes B))  \ar[dd]^= \\
\textup{Tot}(\Gamma_c(U,G(A)) \otimes \Gamma_c(U,G(B))) \ar@<-0.6cm>[d] \ar@<0.6cm>[u]^= \ar@<1cm>[d]^= \ar@<-1cm>[u] & \\
\textup{Tot}(\Gamma(U,G(A)) \otimes \Gamma_c(U,G(B))) \ar[r] & \Gamma_c(U,G(A \otimes B))
}
\]
commutes up to homotopy.
\end{lem}

\dem{}
Using the explicit descriptions above, one needs to prove that the map $\varphi_{r,s}$ from $\left(\prod_{v \notin U} \Gamma(\widehat{K}_v,G_{r-1}(A)) \oplus \Gamma(U, G_r(A))\right) \otimes \left(\prod_{v \notin U} \Gamma(\widehat{K}_v,G_{s-1}(B)) \oplus \Gamma(U,G_s(B))\right)$ to $\prod_{v \notin U} \Gamma(\widehat{K}_v,G_{r+s-1}(A \otimes B)) \oplus \Gamma(U, G_{r+s}(A \otimes B))$ given by 
\[
(a_{r-1}, a_r) \otimes (b_{s-1},b_{s}) \mapsto ((-1)^r \alpha(a_r) \cup b_{s-1} - a_{r-1} \cup \beta(b_s), 0)
\]
is homotopically trivial.
To prove this, consider the maps
\begin{changemargin}{-3cm}{5cm}
{\small 
\[
\left(\prod_{v \notin U} \Gamma(\widehat{K}_v,G_{r-1}(A)) \oplus \Gamma(U, G_r(A))\right) \otimes \left(\prod_{v \notin U} \Gamma(\widehat{K}_v,G_{s-1}(B)) \oplus \Gamma(U,G_s(B))\right	) \xrightarrow{h_{r,s}} \prod_{v \notin U} \Gamma(\widehat{K}_v,G_{r+s-2}(A \otimes B)) \oplus \Gamma(U, G_{r+s-1}(A \otimes B))
\]}
\end{changemargin}
given by $h_{r,s} : (a_{r-1}, a_r) \otimes (b_{s-1}, b_s) \mapsto (0, (-1)^r a_{r-1} \cup b_{s-1})$. Then these maps define an homotopy equivalence between the map $\bigoplus_{r+s = n} \varphi_{r,s}$ and the zero map, proving the lemma.
\enddem

We now prove that the pairing is compatible with coboundary maps in cohomology coming from short exact sequences: 
\begin{lem} \label{lem cup cob}
Let $0 \to A \to B \to C \to 0$ and $0 \to C' \to B' \to A' \to 0$
be two exact sequences of fppf sheaves
on $U$,
and let $B \otimes B' \to D$ be a morphism of fppf sheaves. Assume that the induced morphism $A \otimes C' \to D$ is trivial.

Consider the following diagram
\[
\xymatrix{
H^r_c(U,C) \times H^{s+1}(U,C') \ar[r]^-{\cup} \ar@<-1cm>[d]^-{\partial_r} & H^{r+s+1}_c(U,D) \ar[d]^= \\
H^{r+1}_c(U,A) \times H^s(U,A') \ar[r]^-{\cup} \ar@<-1.2cm>[u]_-{\partial'_s} & H^{r+s+1}_c(U,D) \, , 
}
\]
where the horizontal morphisms are induced by the pairings in Lemma \ref{lem cup} and by the morphism $B \otimes B' \to D$, and the vertical maps are the coboundary morphisms.

Then for all $c \in H^r_c(U,C)$ and $a' \in H^s(U,A')$, we have
\[\partial_r(c) \cup a' + (-1)^r c \cup \partial'_s(a') = 0 \, . \]
\end{lem}

\dem{}
For all fppf sheaves $E$, let $\partial_i^E : G_i(E) \to G_{i+1}(E)$ denote the coboundary map in the Godement complex $G(E)$.

Consider the diagram
induced by $B \otimes B' \to D$:
\[
\xymatrix{
\Gamma(U,G_r(B)) \otimes \Gamma(U,G_{s+1}(B')) \ar[r]^-{\cup} \ar@<-1.4cm>[d]^{\partial^{B}_r} & \Gamma(U, G_{r+s+1}(D)) \ar[d]^= \\
\Gamma(U,G_{r+1}(B)) \otimes \Gamma(U,G_{s}(B')) \ar[r]^-{\cup} \ar@<-1.4cm>[u]_-{\partial^{B'}_s} & \Gamma(U, G_{r+s+1}(D)) \, ,
}
\]
together with the similar diagrams over $\spec \widehat{K}_v$, for all $v \notin S$.

By compatibility of the Godement resolution with tensor product 
(cf.~\cite{friedsus}, Appendix A),
the pairing $\textup{Tot}(G(B) \otimes G(B')) \to G(D)$ is a morphism of complexes. Hence for all $b \in \Gamma(U,G_r(B))$ and $b' \in \Gamma(U,G_{s}(B'))$, we have
\[\partial^{B}_r(b) \cup b' + (-1)^r b \cup \partial^{B'}_s(b') = \partial^D_{r+s}(b \cup b') \, .\]
This formula, its analogue over $\spec \widehat{K}_v$ for $v \notin S$, together with the definition of the
connecting
maps in cohomology via Godement resolutions (recall that for all $n$, the functor $\mathcal{F} \mapsto G_n(\mathcal{F})$ is exact, see \cite{SGA4}, XVII, Proposition 4.2.3), implies Lemma \ref{lem cup cob}.
\enddem

\begin{lem} \label{continuous cup}
Let $N$ be a finite flat commutative 
$U$-group scheme of order $n$, then the pairings \eqref{cup}
$$H^r_c(U,N) \times H^s(U,N^D) \to H^{r+s}_c(U, \mu_n)$$
are continuous.
\end{lem}

\dem{}
The pairings \eqref{cup} are defined via the cup-product on $U$ and the local duality pairings 
$H^{a}(\widehat{K}_v,N) \times H^{b}(\widehat{K}_v, N^D) 
\to H^{a+b}(\widehat{K}_v, \mu_n)$. 
These local pairings are continuous
(see \cite{Ces}, Theorems 5.11 and 6.5). Hence the lemma follows from the definition of the topologies on the cohomology
groups (see section~\ref{section topo}).
\enddem

\begin{rema}
{\rm In \cite{MilADT} (see for example Th.~III.3.1), 
the pairings are defined via the $\ext$ groups, which is quite convenient for 
the definition itself but makes the required commutativities 
of diagrams more difficult to check. Nevertheless, Proposition V.1.20 in \cite{MilEC} provides a similar comparison between both definitions: see the details in the Appendix, Proposition~\ref{compext}.
}
\end{rema}

In order to prove Theorem \ref{thm AV}, we now need to show that the induced map $H^{3-r}_c(U,N) \to 
H^{r}(U,N^D)^*$ is an isomorphism (of topological groups) for every finite flat commutative group
scheme $N$ over $U$ and every $r \in \{ 0,1,2,3 \}$ (recall that the groups 
$H^{r}(U,N^D)$ are equipped with the discrete topology).

\smallskip

We first recall the following lemma (\cite{MilADT}, Lemma III.8.3):
\begin{lem}\label{lem devissage}
Let 
\[0 \to N' \to N \to N'' \to 0 \]
be an exact sequence of finite flat commutative
group schemes over $U$.
If Theorem \ref{thm AV} is true for $N'$ and $N''$, then it is true for $N$.
\end{lem}

\dem{}
Using Proposition \ref{prop supp compact}, 2., the exactness of Pontryagin duality for discrete groups and the pairing in Lemma \ref{lem cup}, one gets
a diagram
of long exact sequences:
\[
\xymatrix{
\dots \ar@{}|(.5){\star}[dr] \ar[r] & H^{3-r}_c(U,N') \ar[r] \ar[d] & H^{3-r}_c(U,N) \ar[r] \ar[d] & H^{3-r}_c(U,N'') \ar[r] \ar[d] \ar@{}|(.5){\star}[dr] & \dots \\
\dots \ar[r] & H^r(U, {N'}^D)^* \ar[r] & H^r(U, {N}^D)^* \ar[r] & H^r(U, {N''}^D)^* \ar[r] & \dots \, .
}
\]
The functoriality of the pairing (see Lemma \ref{lem cup}) implies that both central squares in the diagram are commutative. Lemma \ref{lem cup cob} implies that both extreme squares (denoted~$\star$) are commutative up to sign. By
Lemma~\ref{strictprof} 1., Prop.~\ref{rem topo} and Lemma \ref{continuous cup}, all the maps in this diagram are continuous. Hence the lemma follows
from the five-lemma. 
\enddem

We now want to show that it is sufficient to prove Theorem~\ref{thm AV} for a smaller open subset $V \subset U$. To do this, we need to check the compatibility of the pairing in Theorem \ref{thm AV} with restriction to an open subset of $U$ and with the local duality pairing (see Lemma \ref{lem commutative} below). 

We first define the maps that appear in this lemma. Let $\mathcal{F}$ be a flat affine commutative $U$-group scheme of finite type and let $V \subset U$ be a non empty open subset. Let $W$ denote $U \setminus V$. In diagram \eqref{big diag} below, the first column is the long exact sequence of Proposition \ref{prop supp compact}, 3., and the second column is the localization exact sequence from Prop.~\ref{bonus}.  
The horizontal pairings are either the local duality pairings from \cite{MilADT}, Theorem III.7.1 (first and last rows), using the same sign convention as in the pairing \eqref{cup complexes right}, or the global pairings from Lemma \ref{lem cup} (second and third rows). The proof of Proposition \ref{bonus} provides an isomorphism $H^3_W(U,\G) \cong \bigoplus_{v \in W} H^3_v(\widehat{\mathcal{O}_v}, \G)$, and the natural morphism of complexes $\Gamma_W(U,I^\bullet(\G)) \to \Gamma_c(V,I^\bullet(\G))$ gives a morphism $H^3_W(U,\G) \to H^3_c(V, \G)$, whence natural morphisms $\bigoplus_{v \in W} H^3_v(\widehat{\mathcal{O}_v}, \G) \to H^3_c(V, \G)  \to H^3_c(U, \G)$.


\begin{lem} \label{lem commutative}
Let $\mathcal{F}$, $\mathcal{G}$ be flat affine commutative group schemes of finite type on $U$, together with a pairing $\mathcal{F} \otimes \mathcal{G} \to \G$. Let $V \subset U$ be a non empty open subscheme and $W := U \setminus V$. Then the following diagram
\begin{equation} \label{big diag}
\begin{gathered}
\xymatrix{
\bigoplus_{v \in W} H^{2-r}(\widehat{\mathcal{O}_v}, \mathcal{F}) \times \bigoplus_{v \in W} H^{r+1}_v(\widehat{\mathcal{O}_v}, 
\mathcal{G})  \ar@<-1.2cm>[d] \ar[r] & \bigoplus_{v \in W} H^3_v(\widehat{\mathcal{O}_v}, \G) \ar[d] \\ 
H^{3-r}_c(V, \mathcal{F}) \times H^{r}(V, \mathcal{G}) \ar[r] \ar@<-1.2cm>[d] \ar@<-1.2cm>[u]  & H^3_c(V, \G) \ar[d]^\sim \\ 
H^{3-r}_c(U, \mathcal{F}) \times H^{r}(U, \mathcal{G}) \ar[r] \ar@<-1.2cm>[u] \ar@<-1.2cm>[d] & H^3_c(U, \G) \\
\bigoplus_{v \in W} H^{3-r}(\widehat{\mathcal{O}_v}, \mathcal{F}) \times \bigoplus_{v \in W} H^{r}_v(\widehat{\mathcal{O}_v},  \mathcal{G})  \ar@<-1.2cm>[u] \ar[r] & \bigoplus_{v \in W} H^3_v(\widehat{\mathcal{O}_v}, \G) \ar[u] 
}
\end{gathered}
\end{equation}
is commutative.

In addition, if ${\cal F}$ and $\mathcal{G}$ are finite and flat group schemes, then all the maps in the diagram are continuous.
\end{lem}

\dem{}
\begin{enumerate}
	\item 
We first prove the commutativity of the top rectangle. It is sufficient to prove that the following diagram commutes:
\[
\xymatrix{
\textup{Tot}\left(\bigoplus_{v \in W} \Gamma(\widehat{\mathcal{O}_v},G(\mathcal{F}))[-1] \otimes \bigoplus_{v \in W} \Gamma_v(\widehat{\mathcal{O}_v},G(\mathcal{G}))[1]\right) \ar[r] \ar@<-1cm>[d] & \bigoplus_{v \in W} \Gamma_v(\widehat{\mathcal{O}_v}, G(\G)) \\
\textup{Tot}\left(\bigoplus_{v \in W} \Gamma(\widehat{K}_v, G(\mathcal{F}))[-1] \otimes \bigoplus_{v \in W} \Gamma(\widehat{K}_v, G(\mathcal{G}))\right) \ar[r] \ar@<-1cm>[d] \ar@<-1cm>[u] & \bigoplus_{v \in W} \Gamma(\widehat{K}_v, G(\G))[-1] \ar[u] \ar[d] \\
\textup{Tot}\left(\Gamma_c(V,G(\mathcal{F})) \otimes \Gamma(V,G(\mathcal{G}))\right) \ar[r] \ar@<-1cm>[u] & \Gamma_c(V, G(\G)) \\
}
\]
where the vertical maps are the natural ones and the horizontal pairings are defined earlier. The top rectangle is commutative because of the definition of the pairing involving cohomology with support in a closed subscheme, taking into account the sign conventions in Proposition~\ref{tensorcomplex} in the Appendix. The bottom one is commutative by definition of the pairing involving compact support cohomology.

Assume now that $\mathcal{F}$ and $\mathcal{G}$ are finite flat group schemes. 
Then the following maps are continuous: the 
pairing $H^2_c(V, \mathcal{F}) \times H^1(V, \mathcal{G}) \to H^3_c(V, \Gm)$
(see Lemma \ref{continuous cup}), the pairing
 $H^1(\widehat{\mathcal{O}_v}, \mathcal{F}) \times H^2_v(\widehat{\mathcal{O}_v}, \mathcal{G}) \to H^3_v(\widehat{\mathcal{O}_v}, \Gm)$
(see \cite{MilADT}, Theorem III.7.1) and the map $\bigoplus_{v \in W} H^1(\widehat{\mathcal{O}_v}, \mathcal{F}) \to H^2_c(V, \mathcal{F})$ (see Proposition \ref{rem topo}).

	\item We now prove the commutativity of the rectangle in the middle. Let $\widetilde{\Gamma}(U,G(\mathcal{F})) := \textup{Cone}(\Gamma(U,G(\mathcal{F})) \to \bigoplus_{v \notin U} \Gamma(\widehat{K_v},G(\mathcal{F})) \oplus \bigoplus_{v \in U \setminus V} \Gamma(\widehat{\mathcal{O}_v}, G(\mathcal{F}))) [-1]$. Then functoriality of the cone gives a commutative 
diagram (similar to (\ref{bigdiag}), where $I^{\bullet}(\mathcal F)$ is
replaced by $G(\mathcal{F})$ and by $G(\G)$) of complexes of abelian groups:
\[
\xymatrix{
\textup{Tot}(\Gamma_c(V,G(\mathcal{F})) \otimes \Gamma(V,G(\mathcal{G}))) \ar[r]  & \Gamma_c(V,G(\G)) \\
\textup{Tot}(\widetilde{\Gamma}(U,G(\mathcal{F})) \otimes \Gamma(U,G(\mathcal{G}))) \ar[r] \ar@<-1.2cm>[u] \ar@<1.2cm>[d]^= \ar@<0.8cm>[u]^q  \ar@<-0.8cm>[d] & \widetilde{\Gamma}(U,G(\G)) \ar[d] \ar[u]^q \\
\textup{Tot}(\Gamma_c(U,G(\mathcal{F})) \otimes \Gamma(U,G(\mathcal{G}))) \ar[r] & \Gamma_c(U,G(\G)) \, .
}
\]
Here the maps denoted by $q$ are quasi-isomorphisms (see Remark \ref{rem acyclic} and the proof of the third point in Proposition \ref{prop supp compact}, 
which uses Lemma~\ref{lem dense H1}).
This diagram gives a commutative diagram in the derived category of abelian groups (where all the maps are either the natural ones or the ones constructed above):
\[
\xymatrix{
\Gamma_c(V,G(\mathcal{F})) \otimes^{\L} \Gamma(V,G(\mathcal{G})) \ar[r] \ar@<-1cm>[d] & \Gamma_c(V,G(\G)) \ar[d] \\
\Gamma_c(U,G(\mathcal{F})) \otimes^{\L} \Gamma(U,G(\mathcal{G})) \ar[r] \ar@<-1cm>[u] & \Gamma_c(U,G(\G)) \, .
}	
\]
Taking cohomology of this diagram gives a commutative diagram of abelian groups:
\[
\xymatrix{
H^r_c(V,\mathcal{F}) \times H^s(V,\mathcal{G}) \ar[r] \ar@<-1cm>[d] & H^{r+s}_c(V,\G) \ar[d] \\
H^r_c(U,\mathcal{F}) \times H^s(U,\mathcal{G}) \ar[r] \ar@<-1cm>[u] & H^{r+s}_c(U,\G) \, ,
}
\]
which implies the required commutativity. 

The continuity of the maps in the case where $\mathcal{F}$ and $\mathcal{G}$ are finite flat group schemes is a consequence of Lemma \ref{continuous cup} and of Lemma~\ref{continuuv}.

	\item We now need to prove the commutativity of the bottom rectangle. By Lemma \ref{lem exchange left right}, the following diagram
\[
\xymatrix{
\Gamma_c(U,G(\mathcal{F})) \otimes^{\L} \Gamma(U,G(\mathcal{G})) \ar[r] \ar@<-1cm>[d] & \Gamma_c(U,G(\G)) \ar[d]^= \\
\Gamma(U,G(\mathcal{F})) \otimes^{\L} \Gamma_c(U,G(\mathcal{G})) \ar[r] \ar@<-0.9cm>[u] & \Gamma_c(U,G(\G))  
}
\]
commutes in the derived category. Computing cohomology gives a commutative diagram of abelian groups:
\[
\xymatrix{
H^r_c(U,\mathcal{F}) \times H^s(U,\mathcal{G}) \ar[r] \ar@<-1cm>[d] & H^{r+s}_c(U,\G) \ar[d]^= \\
H^r(U,\mathcal{F}) \times H^s_c(U,\mathcal{G}) \ar[r] \ar@<-1cm>[u] & H^{r+s}_c(U,\G)  \, .
}
\]
	
Let $\Gamma_W(U,G(\mathcal{G})) := \textup{Cone}(\Gamma(U,G(\mathcal{G})) \to \Gamma(V,G(\mathcal{G})))[-1]$. In order to prove the required commutativity, it is enough to prove that the natural diagram
\[
\xymatrix{
\Gamma(U,G(\mathcal{F})) \otimes^{\L} \Gamma_c(U,G(\mathcal{G})) \ar[r]  & \Gamma_c(U,G(\G)) \\
\Gamma(U,G(\mathcal{F})) \otimes^{\L} \Gamma_W(U,G(\mathcal{G})) \ar[r] \ar@<0.9cm>[u]^= \ar@<-1cm>[u] & \Gamma_W(U,G(\G)) \ar[u]
}
\]
commutes in the derived category, where the pairing on the bottom row is defined in a similar way as the pairing \eqref{cup complexes right}. Consider the following diagram in the category of complexes:
\begin{changemargin}{-4cm}{5cm}
\[
\xymatrix{
\textup{Tot}(\Gamma(U,G(\mathcal{F})) \otimes \Gamma(U,G(\mathcal{G}))) \ar[rr] \ar[dd] \ar[rd] & & \Gamma(U,G(\G)) \ar[dd] \ar[rd] & \\
& \textup{Tot}(\Gamma(U,G(\mathcal{F})) \otimes \Gamma(V,G(\mathcal{G}))) \ar[dl] \ar[rr]|!{[r];[r]}\hole & & \Gamma(V,G(\G)) \ar[ld] \\
\prod_{v \notin U} \textup{Tot}(\Gamma(U,G(\mathcal{F})) \otimes \Gamma(\widehat{K}_v,G(\mathcal{G}))) \ar[rr] &  & \prod_{v \notin U} \Gamma(\widehat{K_v}, G(\G)) \, . &
}
\]
\end{changemargin}
This diagram is commutative, hence, using Proposition \ref{tensorcomplex}, it induces a commutative
diagram of complexes at the level of cones:
\begin{changemargin}{-4cm}{5cm}
\[
\xymatrix{
& \textup{Tot}(\Gamma(U,G(\mathcal{F})) \otimes \Gamma_c(U,G(\mathcal{G}))) \ar[rr] \ar[dd]|!{[d];[d]}\hole & & \Gamma_c(U,G(\G)) \ar[dd] \\
\textup{Tot}(\Gamma(U,G(\mathcal{F})) \otimes \Gamma_W(U,G(\mathcal{G}))) \ar[rr]  \ar[ru] \ar[rd] & & \Gamma_W(U,G(\G)) \ar[rd] \ar[ru] & \\
& \textup{Tot}(\Gamma(U,G(\mathcal{F})) \otimes \Gamma(U,G(\mathcal{G}))) \ar[rr]  & & \Gamma(U,G(\G)) \, .
}
\]
\end{changemargin}
The commutativity of the upper face of this last diagram concludes the proof.

Assume now that $\mathcal{F}$ and $\mathcal{G}$ are finite flat group schemes. The only possibly non discrete groups in the diagram are $H^2_c(U, \mathcal{F})$ (in the case $r=1)$ and $H^1(\widehat{\mathcal{O}_v}, \mathcal{F})$ (in the case $r=2$). If $r=1$, the pairing $H^2_c(U, \mathcal{F}) \times H^1(U, \mathcal{G}) \to H^3_c(U, \Gm)$ is continuous by Lemma \ref{continuous cup} and $H^2(\widehat{\mathcal{O}_v}, \mathcal{F}) = 0$ for all $v \in W$ (see for instance \cite{MilADT}, Lemma 1.1), hence all maps are continuous in this case. If $r=2$, then the local pairings $H^1(\widehat{\mathcal{O}_v}, \mathcal{F}) \times H^2_v(\widehat{\mathcal{O}_v}, \mathcal{G}) \to H^3_v(\widehat{\mathcal{O}_v}, \Gm)$ are continuous by \cite{MilADT}, Theorem III.7.1. All the other maps are obviously continuous.
\end{enumerate}
This finishes the proof of Lemma~\ref{lem commutative}.
\enddem

We can now prove the following lemma (see \cite{MilADT}, Lemma III.8.4):
\begin{lem} \label{lem restr}
Let $V \subset U$ be a non empty open subscheme. Let $N$ be a finite flat
commutative 
group scheme over $U$. Then Theorem \ref{thm AV} holds for $N$ over $U$ if and only if it holds for $N_{\vert V}$ over $V$.
\end{lem}

\dem{}
Proposition \ref{prop supp compact},3., Proposition \ref{bonus}, 
Proposition~\ref{rem topo} 
and Lemma \ref{lem commutative} give a commutative diagram of long exact sequences of topological groups:
\[
\xymatrix{
\dots \ar[r] & H^{3-r}_c(V,N) \ar[r] \ar[d] & H^{3-r}_c(U,N) \ar[r] \ar[d] & \bigoplus_{v \in U \setminus V} H^{3-r}(\widehat{\mathcal{O}_v}, N) \ar[r] \ar[d] & \dots \\
\dots \ar[r] & H^r(V, N^D)^* \ar[r] & H^r(U, N^D)^* \ar[r] &  \bigoplus_{v \in U \setminus V} H^r_v(\widehat{\mathcal{O}_v}, N^D)^* \ar[r] & \dots \, , 
}
\]
where the vertical maps are defined via the pairings \eqref{cup} and the local duality pairings of \cite{MilADT}, III.7.1. In particular, the maps $H^{3-r}(\widehat{\mathcal{O}_v},N) \to H^r_v(\widehat{\mathcal{O}_v}, N^D)^*$ are isomorphisms by \cite{MilADT}, Theorem III.7.1.
The middle vertical map is strict by Lemma~\ref{continuous cup} and 
Lemma~\ref{strictprof}.
Therefore the five-lemma gives the result.
\enddem

The end of the proof of Theorem \ref{thm AV}
(which implies in particular 
that by duality the groups $H^r(U,N^D)$ are zero for $r \geq 4$, resp.~for $r=3$ if $U \neq X$) 
is exactly the same as the end of the proof of Theorem III.8.2
in \cite{MilADT}: let $U \subset X$ be a non empty open subset and $N$ be a
finite flat commutative
group scheme over $U$.
\begin{itemize}
	\item If the order of $N$ is prime to $p$, then Theorem \ref{thm AV} is a consequence of Proposition \ref{prop supp compact}, 4.~and \'etale Artin--Verdier duality (Corollary II.3.3 in \cite{MilADT} or Theorem 4.6 in \cite{GS}). Note that it requires to compare the pairing defined in Lemma \ref{lem cup} with the Artin--Verdier pairing using Ext groups as defined in \cite{MilADT} or \cite{GS} : this is explained for instance in Proposition~\ref{compext} of the 
Appendix. Hence by Lemma \ref{lem devissage}, it is sufficient to prove Theorem \ref{thm AV} when the order of $N$ is a power of $p$. 
	\item If the order of $N$ is a power of $p$, the proof of Lemma \ref{profiniteh2} implies that $N$ admits a composition series such that the generic fiber of each quotient is either of height one or the dual of a group of height one. By Lemma \ref{lem devissage}, it is therefore sufficient to prove Theorem \ref{thm AV} in the case $N_K$ or $N_K^D$ have height one.
	\item If $N_K$ or $N_K^D$ have height one, Proposition B.4 and Corollary B.5 in \cite{MilADT} imply that there exists a non empty open subset $V \subset U$ such that $N_{\vert V}$ extends to a finite flat commutative 
group scheme $\widetilde{N}$ over the proper $k$-curve $X$, such that $\widetilde{N}$ or $\widetilde{N}^D$ have height one. Using Lemma \ref{lem restr} twice, it is enough to prove Theorem \ref{thm AV} when $U=X$ and $N$ or $N^D$ have height one.
		\item Lemma III.8.5 in \cite{MilADT} proves Theorem \ref{thm AV} for $U = X$ and $N$ (resp.~$N^D$) of height 
one, by reduction to the classical Serre duality for vector bundles over the smooth projective curve $X$. Indeed, Proposition V.1.20 in \cite{MilEC} proves that the pairings $R\Gamma(X,\mathcal{F}^D) \otimes^\L R\Gamma(X,\mathcal{F}) \to R\Gamma(X,\G)$ defined via Godement resolutions in the proof of Lemma \ref{lem cup} are compatible with the classical pairings using Ext groups that appear in Serre duality.
\enddem
\end{itemize}

As observed in \cite{MilADT}, \S III.8 (remark before Lemma~8.9), 
the group $H^1(U,N)$ is in general infinite if $U \neq X$ and by duality, 
the same is true for $H^2_c(U,N)$. However, the situation is better for 
$H^2$ and $H^1_c$~:

\begin{cor} \label{corh2}
Let $N$ be a finite and flat commutative 
group scheme over
a non empty Zariski open subset $U$ of $X$.
Then the groups $H^2(U,N)$ and $H^1_c(U,N)$ are finite. 
\end{cor}

\dem{} The statement about $H^1_c(U,N)$ is Corollary~\ref{finih1c}. 
The finiteness of $H^2(U,N)$ follows by the duality Theorem~\ref{thm AV}.
\enddem

The previous corollary can be refined in some cases: 

\begin{prop} \label{hassepr}
Let $N$ be a finite and flat commutative 
group scheme over a non empty
affine open subset
$U \subset X$, such that the
generic fiber $N_K$ is local. Then $H^1_c(U,N)=0$.
\end{prop}

\dem{} By the valuative criterion of properness, the
restriction map $H^1(U,N) \to H^1(K,N)$ is injective. 
It is sufficient to show that 
if we choose $v \not \in U$, the restriction
map $H^1(K,N) \to H^1(\widehat K_v,N)$ is injective when $N_K$ is local.
Indeed this implies that $D^1(U,N)=0$, hence $H^1_c(U,N)=0$ by
exact sequence (\ref{didef}) because
$H^0(\widehat K_v,N)=0$ for every completion $\widehat K_v$ of $K$. 

\smallskip

We also reduce to showing that for every finite 
subextension $L/K$ of $\widehat K_v/K$, the 
restriction map $r : H^1(K,N) \to H^1(L,N)$ is injective 
(indeed a $K$-torsor under the finite $K$-group 
scheme $N_K$ is of finite type over $K$, hence it 
has a point over an extension $K'$ of $K$ if and only
if it has a point over a finite subextension of $K'$).
To do this, we argue as in \cite{Ces}, 
Lemma~5.7 a). Since by \cite{riben}, section F, Th.~2, 
$L$ is a separable extension of $K$, the $K$-algebra 
$E:=L \otimes_K L$ is reduced. As $N_K$ is finite and connected, the group 
$N(E)$ is trivial. Let $C^1:={\bf R}_{E/K} (N \times_K E)$ (where $\bf R$ 
denotes Weil's restriction of scalars) be the 
scheme of $1$-cochains with respect to $L/K$, we obtain that 
$C^1(K)$ is trivial, which in turn implies that $\ker r$ is trivial 
by \cite{Ces}, \S 5.1.
\enddem

\begin{rema}
{\rm The finiteness of $H^1_c(U,N)$ (Cor.~\ref{finih1c}) 
relies on the finiteness
of $D^1(U,N)$ proven in \cite{cesnalms}, Th.~2.9. An alternative argument 
is actually available. By \cite{MilADT}, Lemma III.8.9., we can assume 
that $U \neq X$, namely that $U$ is affine. 
By loc.~cit., Th.~II.3.1.~and
Prop.~\ref{prop supp compact}, 4., we can also
assume that the order of $N$ is a power of $p$. 
Let $N_K$ be the generic
fiber of $N$, it is a finite group scheme over $K$.
By \cite{demgab}, IV, \S 3.5, and 
Prop.~\ref{prop supp compact}, 2., it is
sufficient to prove the required finiteness in the following cases~:
$N_K$ is \'etale, $N_K$ is local with \'etale dual, $N_K=\alpha_p$.
The last two cases are taken care of by Prop.~\ref{hassepr}, so we can 
suppose that $N_K$ is \'etale. Let $V \subset U$ be a non empty 
open subset.
By Prop.~\ref{prop supp compact}, we have an exact sequence
$$H^1_c(V,N) \to H^1_c(U,N) \to \bigoplus_{v \in U \setminus V}
H^1(\widehat{\calo_v},N) .$$
Since the generic fiber of $N$ is \'etale, the group
$H^1(\widehat{\calo_v},N)$ is finite by \cite{MilADT}, Rem.~III.7.6.
(this follows from the fact that $H^1(\widehat{\calo_v},N)$ is a compact
subgroup of the discrete group $H^1(K_v,N)$), hence the finiteness 
of $H^1_c(U,N)$ is equivalent to the finiteness of $H^1_c(V,N)$, which 
in turn is equivalent to the finiteness of $D^1(V,N)$. The latter 
holds for $V$ sufficiently small: either apply 
\cite{gonzcrelle}, Lemma 4.3.~(which relies
on an embedding of $N_K$ into an abelian variety) or 
reduce (as in \cite{MilADT}, Lemma~III.8.9.) to the case when 
$N^D$ is of height one. Indeed by loc.~cit., Cor.~III.B.5., the 
assumption that $N^D$ is of height one implies that 
for $V$ sufficiently small, the restriction of $N$ to $V$ extends to
a finite and flat commutative group scheme
$\widetilde N$ over $X$. 
Then the finiteness of $H^1_c(X,\widetilde N)$
implies the finiteness of $H^1_c(V,\widetilde N)=
H^1_c(V,N)$ by Prop~\ref{prop supp compact}, 3., because the groups
$H^0(\widehat{\calo_v},\widetilde N)$ are finite.
}
\end{rema}

\section{The number field case} \label{sect3}

Assume now that $K$ is a number field and set $X=\spec \calo_K$. 
Let $U$ be a non empty Zariski open subset of $X$. Let $n$ be the 
order of the finite and flat commutative 
group scheme $N$. To prove 
Theorem~\ref{thm AV} in this case, one follows exactly the same method
as in \cite{MilADT}, Th.~III.3.1.~and Cor.~III.3.2.~once 
Proposition~\ref{prop supp compact} has been proved. Namely 
Proposition~\ref{prop supp compact}, 4., shows that on $U[1/n]$, 
Theorem~\ref{thm AV} reduces to the \'etale Artin--Verdier Theorem 
(\cite{MilADT}, II.3.3 or Theorem 4.6 in \cite{GS}).
Here we can use 
a definition of the pairings similar to Lemma~\ref{lem cup},
or a definition via 
the Ext pairings as in loc.~cit.~(the two definitions 
coincide, the argument being the same 
as in Proposition~\ref{compext} of the Appendix).
Now Proposition~\ref{prop supp compact}, 3., gives 
a commutative diagram as in the end of the proof of \cite{MilADT}, Th.
III.3.1.~(with completions $\widehat \calo_v$ 
instead of henselizations $\calo_v$). 
Theorem~\ref{thm AV} follows by the five-lemma, using the result 
on $U[1/n]$ and the local duality statement \cite{MilADT}, Th.~III.1.3.

\begin{rema}
{\rm In the number field case, one can as well (as in \cite{MilADT}, 
\S III.3) work from the 
very beginning with henselizations $\calo_v$ and not with completions 
$\widehat \calo_v$ to define cohomology with compact support. Indeed 
the local theorem (loc.~cit., Th.~III.1.3) still holds with henselian 
(not necessarily complete) discrete valuation ring 
with finite residue field when 
the fraction field is of characteristic zero. Hence the only 
issue here is commutativity of diagrams. Nevertheless, we felt that 
it was more convenient to have a uniform statement 
(Proposition~\ref{prop supp compact}) in both characteristic $0$ and 
characteristic $p$ situations.  
}
\end{rema}

\section[A]{Appendix} \label{appsect}
\subsection{Cone and tensor products} \label{appendix a}

\begin{prop} \label{tensorcomplex}
Let $\mathcal{A}$ be the category of fppf sheaves
over a scheme $T$.
Let $A$, $B$ and $C$ be three complexes in $\mathcal{A}$. Let $f : A \to B$ be a morphism of complexes. Then there are commutative diagrams
(where $\otimes$ denotes the total tensor product of complexes) such that 
the vertical maps are isomorphisms of complexes:
\[
\xymatrix{
A \otimes C \ar[r]^{f \otimes 1} \ar[d]^-= & B \otimes C \ar[d]^-= \ar[r]^-{i \otimes 1} & \textup{Cone}(f) \otimes C \ar[d]^-\sim \ar[r]^-{-\pi \otimes 1} & A[1] \otimes C \ar[d]^-\sim \\
A \otimes C \ar[r]^-{f \otimes 1} & B \otimes C \ar[r]^-{i'} & \textup{Cone}(f \otimes 1) \ar[r]^-{-\pi} & (A \otimes C)[1] \, ,
}
\]
where the vertical isomorphisms involve no signs, and 
\[
\xymatrix{
C \otimes A \ar[r]^-{1 \otimes f} \ar[d]^-= & C \otimes B \ar[d]^-= \ar[r]^-{1 \otimes i} & C \otimes \textup{Cone}(f) \ar[d]^-\sim \ar[r]^-{-1 \otimes \pi} & C \otimes A[1] \ar[d]^-\sim \\
C \otimes A \ar[r]^-{1 \otimes f} & C \otimes B \ar[r]^-{i'} & \textup{Cone}(1 \otimes f) \ar[r]^-{-\pi} & (C \otimes A)[1] \, ,
}
\]
where the two last vertical maps involve a sign $(-1)^r$ on the factor $C_r \otimes A_s$.
\end{prop}

\dem{}  In the first diagram, define
the non obvious map $\textup{Cone}(f) \otimes C \to \textup{Cone}(f \otimes 1)$ (resp.~$A[1] \otimes C \to (A \otimes C)[1]$) by the isomorphism $(B_r \oplus A_{r+1}) \otimes C_s \to (B_r \otimes C_s) \oplus (A_{r+1} \otimes C_s)$ (resp.~by the identity of $A_{r+1} \otimes C_s$). In the second diagram, the non obvious map $C \otimes \textup{Cone}(f) \to \textup{Cone}(1 \otimes f)$ (resp.~$C \otimes A[1] \to (C \otimes A)[1]$) is given by the isomorphism $C_r \otimes (B_s \oplus A_{s+1})  \to (C_r \otimes B_s) \oplus (C_r \otimes A_{s+1})$ that maps $c \otimes (b,a)$ to $(c \otimes b, (-1)^r c \otimes a)$ (resp.~by the automorphism of $C_{r} \otimes A_{s+1}$ given by $c \otimes a \mapsto (-1)^r c \otimes a$). The proposition is then straightforward.
\enddem

\subsection{Comparison of two pairings} \label{appendix b}

Let $U$ be a non empty Zariski open subset of a smooth, projective,
geometrically integral curve defined over a finite field.

\begin{prop} \label{compext}
Let $A$, $B$ and $C$ be three fppf sheaves on $U$, endowed with a pairing $A \otimes B \to C$. Then there is a commutative diagram
\[
\xymatrix{
H^r(U,A) \otimes H^s_c(U,B) \ar[r] \ar[d] & H^{r+s}_c(U,C) \ar[d]^= \\
\textup{Ext}^r_U(B,C)\otimes H^s_c(U,B) \ar[r] & H^{r+s}_c(U,C) \, ,
}
\]
where the top pairing is the one from \eqref{cup complexes right} and the bottom one is the pairing from \cite{MilADT}, Proposition III.0.4.e. The same holds 
for \'etale sheaves instead of fppf sheaves if we replace
fppf cohomology (resp.~compact support fppf cohomology) by 
\'etale cohomology (resp.~compact support \'etale cohomology);
in the \'etale case the bottom pairing is the one from 
loc.~cit., Proposition II.2.5.~(or \cite{GS}).
\end{prop}

\dem{} We prove the statement for fppf sheaves (the \'etale case is similar).
Consider the natural morphisms of complexes:
\[
\textup{Tot}(G(A) \otimes G(B)) \to G(A \otimes B) \to G(C) \, .
\]
Using \cite[\href{http://stacks.math.columbia.edu/tag/0A90}{Tag 0A90}]{SP}, one gets a natural morphism of complexes $G(A) \to \mathcal{H}om^\bullet(G(B),G(C))$ and a commutative diagram of complexes:
\[
\xymatrix{
\textup{Tot}(G(A) \otimes G(B)) \ar[r] \ar[d] & G(A \otimes B) \ar[r] & G(C) \ar[d]^= \\
\textup{Tot}(\mathcal{H}om^\bullet(G(B),G(C)) \otimes G(B)) \ar[rr] & & G(C) \, ,
}
\]
where the second pairing is the natural one. All morphisms in this diagram involve no extra-sign.

Let $G(C) \to I$ be an injective resolution. Then one gets a commutative diagram
\[
\xymatrix{
\textup{Tot}(G(A) \otimes G(B)) \ar[r] \ar[d] & G(A \otimes B) \ar[r] & G(C) \ar[d] \\
\textup{Tot}(\mathcal{H}om^\bullet(G(B),I) \otimes G(B)) \ar[rr] & & I \, .
}
\]
Taking global sections, one gets a commutative diagram:
\begin{equation} \label{diag comp cup-Ext}
\begin{gathered}
\xymatrix{
\textup{Tot}(\Gamma(U,G(A)) \otimes \Gamma(U,G(B))) \ar[r] \ar[d] & \Gamma(U,G(C)) \ar[d]^\sim \\
\textup{Tot}(\textup{Hom}_U^\bullet(G(B),I) \otimes \Gamma(U,G(B))) \ar[r] & \Gamma(U,I) \, .
}
\end{gathered}
\end{equation}
Taking cohomology, one gets a commutative diagram comparing the pairing from the beginning of section \ref{sect2} to the classical Ext-pairing:
\[
\xymatrix{
H^r(U,A) \otimes H^s(U,B) \ar[r] \ar[d] & H^{r+s}(U,C) \ar[d]^= \\
\textup{Ext}^r_U(B,C) \otimes H^s(U,B) \ar[r] & H^{r+s}(U,C) \, .
}
\]

Applying functoriality of cone to \eqref{diag comp cup-Ext} and to the similar pairing over completions of $K$, one gets a commutative diagram of complexes
\[
\xymatrix{
\textup{Tot}(\Gamma(U,G(A)) \otimes \Gamma_c(U,G(B))) \ar[r] \ar[d] & \Gamma_c(U,G(C)) \ar[d]^\sim \\
\textup{Tot}(\textup{Hom}_U^\bullet(G(B),I) \otimes \Gamma_c(U,G(B))) \ar[r] & \Gamma_c(U,I) \, .
}
\]
Taking cohomology, we get the required commutative diagram.
\enddem

\begin{rema}  {\rm \label{rem appendix} \hfill
\begin{enumerate}
	\item A similar diagram holds with compact support cohomology groups on the left and Ext-groups on the right. In this case, one gets a commutative diagram
\[
\xymatrix{
H^r_c(U,A) \otimes H^s(U,B) \ar[r] \ar[d] & H^{r+s}_c(U,C) \ar[d]^= \\
H^r_c(U,A) \otimes \textup{Ext}^s_U(A,C) \ar[r] & H^{r+s}_c(U,C) \, ,
}
\]
where the first pairing is the one from Lemma \ref{lem cup}, while the vertical map and the bottom pairing both involve a $(-1)^{r s}$ sign.
	\item Similar commutative diagrams hold over an arbitrary basis, with compact support cohomology replaced by cohomology with support in a closed subscheme (with a similar proof).
\end{enumerate}}
\end{rema}

\merci We thank K.~\v{C}esnavi\v{c}ius, C.~P\'epin, A.~Schmidt and T.~Szamuely 
for enlightening comments on the first draft of this paper. We also thank T.~Suzuki for sharing his preprint \cite{suzuki}. 
We warmly thank the anonymous referee for his/her very thorough reading 
of the article and for having made numerous interesting suggestions.

\bigskip

\noindent
Sorbonne Universit\'e, Universit\'e Paris Diderot, CNRS, Institut de Math\'ematiques de Jussieu-Paris Rive Gauche, IMJ-PRG, F-75005, Paris, France\\
and D\'epartement de math\'ematiques et applications, \'Ecole normale sup\'erieure, CNRS, PSL Research University, 45 rue d'Ulm, 75005 Paris, France
\smallskip

\noindent
cyril.demarche@imj-prg.fr

\bigskip

\noindent
Laboratoire de Math\'ematiques d'Orsay, Univ.~Paris-Sud, CNRS,
Universit\'e Paris-Saclay, 91405 Orsay, France.

\smallskip

\noindent
david.harari@math.u-psud.fr

\end{document}